\documentclass[a4paper,11pt,reqno]{amsart}

\usepackage{amsmath,amstext,amsbsy,amsfonts,eucal,latexsym,amssymb,dsfont,stmaryrd}
\usepackage{graphicx}
\usepackage{float}
\usepackage{multicol}
\usepackage{subfigure}
\usepackage{color}
\usepackage{indentfirst}

\newcommand{\VE}{{\mathbf{E}}}

\newcommand{\VR}{{\mathbf{R}}}

\newcommand{\Dsf}{\mathsf{D}}

\providecommand{\FL}{\mathfrak{L}}

\providecommand{\FR}{\mathfrak{R}}

\newcommand{\eps}{{\varepsilon}}
\newcommand{\nnabla}{\mathcal{D}}

\usepackage{comment}

\numberwithin{equation}{section}
\newtheorem{theorem}{Theorem}
\newtheorem{lemma}{Lemma}
\newtheorem{remark}{Remark}

\newtheorem{proposition}{Proposition}
\newtheorem{assumption}{Assumption}

\newcommand{\out}{\text{out}}
\newcommand{\inn}{\text{in}}

\begin{document}
\title{Automated computation of topological derivatives with application to nonlinear elasticity and reaction-diffusion problems}

\author{Peter Gangl}
\address[P. Gangl]{TU Graz, Steyrergasse 30/III, 
    8010 Graz, Austria}
\email{gangl(at)math.tugraz.at}

\author{Kevin Sturm}
\address[K. Sturm]{TU Wien, Wiedner Hauptstr. 8-10,
      1040 Vienna, Austria}
\email{kevin.sturm(at)tuwien.ac.at}

\date{\today}

\begin{abstract}
$\phantom{oi}$ \\
\textbf{Purpose}
While topological derivatives have proven useful in applications of topology optimisation and inverse problems, their mathematically rigorous derivation remains an ongoing research topic, in particular in the context of nonlinear partial differential equation (PDE) constraints.\\
\textbf{Design/methodology/approach}
We present a systematic yet formal approach for the computation of topological derivatives of a large class of PDE-constrained topology optimization problems with respect to arbitrary inclusion shapes. Scalar and vector-valued as well as linear and nonlinear elliptic PDE constraints are considered in two and three space dimensions including a nonlinear elasticity model and nonlinear reaction-diffusion problems. The systematic procedure follows a Lagrangian approach for computing topological derivatives. \\
\textbf{Findings} For problems where the exact formula is known, the numerically computed values show good coincidence. Moreover, by inserting the computed values into the topological asymptotic expansion, we verify that the obtained values satisfy the expected behaviour also for other, previously unknown problems, indicating the correctness of the procedure. \\
\textbf{Originality/value} We present a systematic approach for the computation of topological derivatives that is applicable to a large class of problems. Most notably, our approach covers the topological derivative for a nonlinear elasticity problem, which has not been reported in the literature.
\end{abstract}


\maketitle

\pagestyle{myheadings} \thispagestyle{plain} \markboth{}{}

 \section{Introduction}

The topological derivative concept was first used for finding optimal locations of holes in mechanical structures in \cite{EschenauerKobelevSchumacher1994} and was later introduced in a mathematically concise way in the publications \cite{SokolowskiZochowski1999} and \cite{GarreauGuillaumeMasmoudi2001}. Given a shape function $\mathcal J$ that maps a shape $\Omega$ to a real number $\mathcal J(\Omega)$, the topological derivative at a spatial point $z \in \Omega$ measures the sensitivity of $\mathcal J$ with respect to a small topological perturbation of the shape $\Omega$. Denoting the perturbed shape by $\Omega_\eps$, e.g., $\Omega_\eps := \Omega \setminus \omega_\eps$ with $\omega_\eps = B_\eps(z)$, the topological derivative is defined as
\begin{align}
    d \mathcal J(\Omega)(z) := \underset{\eps \searrow 0}{\mbox{lim }} \frac{\mathcal J(\Omega_\eps)-\mathcal J(\Omega)}{|\omega_\eps|},
\end{align}
thus satisfying a topologically asymptotic expansion of the form
\begin{align} \label{eq_topAsympExpIntro}
    \mathcal J(\Omega_\eps) = \mathcal J(\Omega) + |\omega_\eps| d\mathcal J(\Omega)(z) + o(\eps) \quad \mbox{as } \eps \searrow 0.
\end{align}

Since its introduction, the topological derivative concept has found application mostly in the context of topology optimization for engineering applications by level set approaches \cite{AmstutzAndrae:2006a, AllaireJouve:2006a, BurgerHacklRing:2004a}, but has also been utilized in medical applications such as electrical impedance tomography (EIT) \cite{HintermuellerLaurainNovotny2012} or mathematical image processing \cite{Hintermueller2005}.

Most practically relevant engineering applications involve a partial differential equation (PDE) constraint and thus are of the type
\begin{align}   \label{eq_intro_pdeconstr}
    \underset{\Omega}{\mbox{inf }} J(\Omega, u) \qquad \mbox{subject to } \; e(\Omega; u) = 0 
\end{align}
with a PDE operator $e(\Omega; \cdot)$. In many cases (in particular when the principal part of the PDE operator is subject to a topological perturbation), the derivation of topological derivatives involves an asymptotic analysis of the state variable of the form
\begin{align} \label{eq_expansion_ueps}
    u_\eps(z+\eps x) = u_0(z+\eps x) + \eps K(x) + o(\eps) \quad \mbox{ as } \eps \searrow 0,
\end{align}
where $u_\eps$ is defined by the perturbed PDE operator, $e(\Omega_\eps, u_\eps) = 0$, and $K$ is a corrector function defined as the solution to a transmission problem on the unbounded domain $\VR^d$. Often, also a similar asymptotic analysis of an adjoint variable is needed.
There exist different approaches for the derivation of topological derivatives for PDE-constrained optimization problems of the type \eqref{eq_intro_pdeconstr}. We mention the approaches by Novotny and Sokolowski \cite{MR3013681}, the approach by Amstutz \cite{Amstutz2006} as well as the averaged adjoint approach \cite{a_ST_2020a}. Finally, we mention the approach that was introduced by Delfour in \cite{c_DE_2018b} and applied to quasilinear problems in \cite{a_GAST_2020a, GanglSturm2021Hcurl}, which does not involve an asymptotic analysis of an adjoint variable.

For many problems involving linear or semilinear PDE constraints, topological derivatives are rather well-understood. In the cases where the shape of the inclusion $\omega$ is a disk or ellipse in 2D or a ball or ellipsoid in 3D, the corrector function $K$ can be determined analytically and the topological derivative can be obtained in a closed form involving so-called polarisation tensors \cite{AmmariKang2007}. When $\omega$ has a more general shape as well as in the case of quasilinear PDE constraints, explicit formulas have not been reported in the literature. Topological derivatives in the context of quasilinear PDE constraints were first treated theoretically in \cite{a_AMBO_2017a}, and later also in \cite{a_GAST_2020a}. The numerical computation for a quasilinear problem was first discussed in \cite{a_AMGA_2019a} in the context of two-dimensional nonlinear magnetostatics, and later in the context of three-dimensional nonlinear magnetostatics in an $H(curl)$ setting in \cite{GanglSturm2021Hcurl}.

In this paper, we make the assumption that an expansion of the form \eqref{eq_expansion_ueps} holds and apply the systematic approach of M.C. Delfour \cite{c_DE_2018b}, see also \cite{a_GAST_2020a}, to a large class of problems involving linear and nonlinear elliptic PDEs as constraints. The suggested procedure covers scalar and vector-valued problems in two and three space dimensions. In particular, we will also treat the problem of nonlinear elasticity in two space dimensions using a St. Venant-Kirchhoff material which, to the best of the authors' knowledge, has not been reported on in the literature so far.

We stress that our generic automated approach is formal since a rigorous derivation would require a detailed asymptotic analysis of the variation of the state \eqref{eq_expansion_ueps}. Usually, this entails to show a sufficiently fast decay of the corrector $K$ as $|x| \rightarrow \infty$. While the computation is formal, we show numerically that, for the considered problems, the computed values satisfy a topologically asymptotic expansion of the form \eqref{eq_topAsympExpIntro}.

 \section{Lagrangian approach for computing topological derivatives} \label{sec_derivTD}
 In this section, we present the general approach for computing topological derivatives of a class of model problems based on the Lagrangian approach introduced in \cite{c_DE_2018b}. The considered problems include linear and nonlinear, scalar and vector-valued elliptic PDE constraints in two or three space dimensions.
 
 \subsection{Class of considered problems}
We consider a PDE-constrained shape optimization problem in $d$ space dimensions where the solution to the PDE is $\VR^m$-valued. This covers the case of scalar quantities, $m=1$ or vector-valued problems such as elasticity, $m=d$. For a function $\varphi : \VR^d \rightarrow \VR^m$, we denote by $\nnabla \varphi \in \VR^{m \times d}$ its Jacobian,  $(\nnabla \varphi)_{i,k} = \frac{\partial \varphi_i }{\partial x_k}$ for $i \in \{1, \dots m\}, k \in \{1,\dots d\}$. We recall the Euclidean vector product $a \cdot b = \sum_{i=1}^m a_i b_i$ for $a, b \in \VR^m$ and will also use this notation for $m=1$, noting that the vector product is just a simple product then. Moreover, we denote by $A:B = \sum_{i=1}^m\sum_{k=1}^d A_{i,k} B_{i,k}$ the Frobenius inner product of two matrices $A, B \in \VR^{m \times d}$ and again note that the Frobenius inner product reduces to the Euclidean inner product of two vectors when $m=1$.

We consider a computational domain $\Dsf$ which is subdivided into two open disjoint subdomains, $\Dsf = \Omega^{\inn} \cup \Omega^{\out}$. We introduce the abbreviation $\Omega := \Omega^{\inn}$, such that $\Omega^{\out} = \Dsf \setminus \overline \Omega$. Moreover, we consider four operators
\begin{align*}
    A_1^{\inn}, A_1^{\out}:&\, \Dsf \times \VR^m \times \VR^{m \times d} \rightarrow \VR^m \\
    A_2^{\inn}, A_2^{\out}:&\, \Dsf \times \VR^m \times \VR^{m \times d} \rightarrow \VR^{m \times d}
\end{align*}
to define the two piecewise defined operators $A_1^{\Omega}:\Dsf \times \VR^m \times \VR^{m \times d} \rightarrow \VR^m$, $A_2^{\Omega}:\Dsf \times \VR^m \times \VR^{m \times d} \rightarrow \VR^{m \times d}$,
\begin{align}
    &A_1^{\Omega}(x, y_1, y_2) := \chi_\Omega(x) A_1^{\inn}(x, y_1, y_2) + \chi_{\Dsf \setminus \Omega}(x) A_1^{\out}(x, y_1, y_2),\\
    &A_2^{\Omega}(x, y_1, y_2) := \chi_\Omega(x) A_2^{\inn}(x, y_1, y_2) + \chi_{\Dsf \setminus \Omega}(x) A_2^{\out}(x, y_1, y_2),
\end{align}
which will represent the left hand side of an abstract PDE constraint. The right hand side will comprise $F_1^\Omega : \Dsf \rightarrow \VR^m$, $F_2^\Omega : \Dsf \rightarrow \VR^{m\times d}$ with
\begin{align}
    &F_i^{\Omega}(x) := \chi_\Omega(x) F_i^{\inn}(x) + \chi_{\Dsf \setminus \Omega}(x) F_i^{\out}(x),    
\end{align}
for $i =1, 2$ with functions $F_1^{\inn}, F_1^{\out} : \Dsf \rightarrow \VR^m$ and $F_2^{\inn}, F_2^{\out} :\Dsf \rightarrow \VR^{m \times d}$. Finally, we assume that the boundary $\partial \Dsf$ of $\Dsf$ is subdivided into two subsets $\Gamma_D, \Gamma_N$ and consider a function $g_N : \Gamma_N \rightarrow \VR^{m}$ to represent inhomogeneous Neumann boundary conditions. For sake of more compact presentation, we only consider homogeneous Dirichlet conditions on the Dirichlet boundary $\Gamma_D$ and remark that an extension to inhomogenous conditions can be obtained by minor modifications.
For a given admissible subdomain $\Omega \subset \Dsf$, we consider the PDE constraint to find $u \in \VE(\Dsf)$ such that
\begin{align} \label{eq_abstr_pde}
    \int_\Dsf  A_1^\Omega(x,u, \nnabla u) \cdot \psi + A_2^\Omega(x,u, \nnabla u) : \nnabla \psi \; \mbox dx = \int_\Dsf F_1^\Omega(x) \cdot \psi  + F_2^\Omega(x) : \nnabla \psi \; \mbox dx + \int_{\Gamma_N} g_N \cdot \psi \;  \mbox dS
\end{align}
for all $\psi \in \VE(\Dsf)$, where $\VE(\Dsf)$ is the function space on which the PDE is posed and which includes the homogeneous Dirichlet conditions on $\Gamma_D$.
Similarly, we consider two functions $j^{\inn}, j^{\out}: \Dsf \times \VR^m \times \VR^{m \times d} \rightarrow \VR$ and, for a given admissible set $\Omega \subset \Dsf$, define the piecewise defined function $j^\Omega: \Dsf \times \VR^m \times \VR^{m \times d} \rightarrow \VR$,
\begin{align}
    j^\Omega(x, y_1, y_2) := \chi_\Omega(x) j^{\inn}(x, y_1, y_2) + \chi_{\Dsf \setminus \Omega}(x) j^{\out}(x, y_1, y_2).    
\end{align}
For $\varphi \in \VE(\Dsf)$, we will consider the cost function 
\begin{align} \label{eq_abstr_J}
    J(\Omega, \varphi, \nnabla \varphi) := J^{\text{vol}}(\Omega,\varphi, \nnabla \varphi)  + J^{\text{bnd}}(\varphi, \nnabla \varphi)
\end{align}
with 
\begin{align}
    J^{\text{vol}}(\Omega, \varphi, \nnabla \varphi) := \int_\Dsf  j^\Omega(x, \varphi, \nnabla \varphi)\; dx,  \qquad
    J^{\text{bnd}}(\varphi, \nnabla \varphi) := \int_{\partial \Dsf}j^{bnd}(x,\varphi,  \nnabla \varphi ) \; dS
\end{align}
for a function $j^{bnd}:\partial \Dsf \times \VR^m \times \VR^{m\times d} \rightarrow \VR$.
Summarizing, we consider the abstract class of PDE-constrained topology optimization problem which can be written as
\begin{align} \label{eq_abstrProblem}
    \begin{aligned}
    &\underset{\Omega \in \mathcal A}{\mbox{min }} J^{\text{vol}}(\Omega,u, \nnabla u)  + J^{\text{bnd}}(u, \nnabla u) \\
    &\mbox{subject to }u  \in \VE(\Dsf) \mbox{ solves } \eqref{eq_abstr_pde}
    \end{aligned}
\end{align}
where $\mathcal A$ denotes the set of admissible subsets of $\Dsf$.

\begin{remark}
    Although we consider only homogeneous boundary conditions, the formulas for inhomogeneous Dirichlet boundary conditions will not affect the formula of the first topological derivative; see, e.g., \cite{a_BAST_2021a} where mixed (Dirichlet and Neumann) inhomogeneous boundary conditions are considered. 
\end{remark}
\subsection{Domain perturbation}
Let $\omega \subset \VR^d$ with $0 \in \omega$ represent the shape of the considered topological perturbation. For a point $z \in \Dsf \setminus \partial \Omega$ and a given small parameter $\eps$, we define $\omega_\eps(z):= z + \eps \omega$ as well as the perturbed domain
\begin{align*}
    \Omega_\eps(z) := \begin{cases}
                        \Omega \setminus \overline \omega_\eps(z), & z \in \Omega, \\
                        \Omega \cup \omega_\eps(z), & z \in \Dsf \setminus \overline \Omega.
                   \end{cases}
\end{align*}
From now on, we assume that a fixed point $z \in \Dsf \setminus \overline \Omega$ is given and 
set $\omega_\eps:= \omega_\eps(z)$, $\Omega_\eps:= \Omega_\eps(z)$. Moreover, we use the abbreviation $\Dsf_\eps := T_\eps^{-1}(\Dsf)$, where $T_\eps(x):= z+ \eps x$.  We also introduce the abbreviating notations $A_1^{(\eps)} := A_1^{\Omega_\eps}$, $A_2^{(\eps)} := A_2^{\Omega_\eps}$, $F_1^{(\eps)} := F_1^{\Omega_\eps}$, $F_2^{(\eps)} := F_2^{\Omega_\eps}$, $j^{(\eps)} := j^{\Omega_\eps}$ and define the perturbed Lagrangian of the abstract optimization problem \eqref{eq_abstrProblem}
\begin{align*}
    G(\eps, \varphi, \nnabla \varphi, \psi, \nnabla\psi) =& \, J(\Omega_\eps, \varphi, \nnabla \varphi) + \langle A_1^{(\eps)}(\varphi, \nnabla\varphi), \psi \rangle + \langle A_2^{(\eps)}(\varphi, \nnabla\varphi), \nnabla\psi \rangle \\
                                                          &- \langle F_1^{(\eps)}, \psi \rangle - \langle F_2^{(\eps)}, \nnabla \psi \rangle -\langle g_N, \psi \rangle,
\end{align*}
where $\langle \cdot, \cdot \rangle$ denotes the duality product in the corresponding spaces.

\subsection{Formal variation of the state variable}
The following derivation is motivated by the rigorous results obtained in \cite{a_GAST_2020a}, where a simplified approach 
to the computation of quasi-linear derivatives is proposed. Here we generalise these results formally to the class of non-linear
problems introduced in \eqref{eq_abstrProblem}. We also refer to \cite{a_AMGA_2019a,a_AMBO_2017a} for another approach to the computation of topological derivatives for quasi-linear problems.

Let $\psi\in \VE(\Dsf)$. The perturbed state equation reads: $u_\eps$ solves
\begin{align} \label{eq_pertState}
    \begin{split}
    \int_\Dsf  A_1^{(\eps)}(x,u_\eps, \nnabla u_\eps) \cdot \psi + A_2^{(\eps)}(x,u_\eps, \nnabla u_\eps) : \nnabla \psi \; dx = \int_\Dsf F_1^{(\eps)}(x) \cdot \psi  &  + F_2^{(\eps)}(x) : \nnabla \psi \; dx \\
                                                                                                                                                                 &  + \int_{\Gamma_N} g_N \cdot \psi \; dS.
\end{split}
\end{align}
Subtracting the state equation \eqref{eq_pertState} for $\eps >0$ from the state equation for $\eps =0$ leads to 
\begin{align*}
    \begin{split}
        \int_\Dsf[  A_1^{(\eps)}(x,u_\eps, \nnabla u_\eps)  - & A_1^{(0)}(x,u_0, \nnabla u_0)] \cdot \psi + [A_2^{(\eps)}(x,u_\eps, \nnabla u_\eps)  - A_2^{(0)}(x,u_0, \nnabla u_0)]  : \nnabla \psi \;dx  \\
                                                                                                     &= \int_{\omega_\eps} (F_1^{\inn}-F_1^{\out})(x) \cdot \psi  + (F_2^{\inn} - F_2^{\out})(x) : \nnabla \psi\; dx.
\end{split}
\end{align*}
Using the change of variables $T_\eps(x)=z+\eps x$, we obtain with the notation $\Dsf_\eps = T_\eps^{-1}(\Dsf)$: 
\begin{align*}
    \int_{\Dsf_\eps} & [  A_1^{(\eps)}(T_\eps(x),u_\eps \circ T_\eps, (\nnabla u_\eps)\circ T_\eps)  - A_1^{(0)}(T_\eps(x),u_0\circ T_\eps, (\nnabla u_0)\circ T_\eps)] \cdot (\psi\circ T_\eps) \\
    +& [A_2^{(\eps)}(T_\eps (x) ,u_\eps\circ T_\eps, (\nnabla u_\eps)\circ T_\eps)  - A_2^{(0)}(T_\eps(x),u_0\circ T_\eps, (\nnabla u_0)\circ T_\eps)] : (\nnabla \psi)\circ T_\eps  \; dx \\
     &= \int_{\omega} ((F_1^{\inn}-F_1^{\out})\circ T_\eps) \cdot  (\psi\circ T_\eps)  + ((F_2^{\inn} - F_2^{\out})\circ T_\eps ) : ((\nnabla \psi )\circ T_\eps)\; dx.
\end{align*}
Using $(\nnabla \varphi)\circ T_\eps = \frac{1}{\eps} \nnabla(\varphi \circ T_\eps)$ and the first variation of the state $K_\eps := \frac{(u_\eps - u_0)\circ T_\eps}{\eps}$, after multiplication with $\eps$, we obtain
\begin{align*}
    \int_{\Dsf_\eps} & [  A_1^{(\eps)}(T_\eps(x),u_0 \circ T_\eps + \eps K_\eps, (\nnabla u_0)\circ T_\eps + \nnabla K_\eps)   \\
                     & \hspace{4cm}- A_1^{(0)}(T_\eps(x),u_0\circ T_\eps, (\nnabla u_0)\circ T_\eps)] \cdot (\eps \psi\circ T_\eps) \\
    +& \int_{\Dsf_\eps}[  A_2^{(\eps)}(T_\eps(x),u_0 \circ T_\eps + \eps K_\eps, (\nnabla u_0)\circ T_\eps + \nnabla K_\eps) \\
     & \hspace{5cm} - A_2^{(0)}(T_\eps(x),u_0\circ T_\eps, (\nnabla u_0)\circ T_\eps)] : \nnabla (\psi\circ T_\eps)  \; dx\\
     &= \int_{\omega} ((F_1^{\inn}-F_1^{\out}) \circ T_\eps) \cdot  (\eps \psi\circ T_\eps)  + ((F_2^{\inn} - F_2^{\out})\circ T_\eps ) : \nnabla (\psi \circ T_\eps) \; dx.
\end{align*}
Now we make the following assumption:
\begin{assumption}\label{ass:main}
    \begin{itemize}
        \item[(i)] We assume that $u_0$ is continuously differentiable at $z$.
        \item[(ii)] Let $K_\eps := \frac{(u_\eps - u_0)\circ T_\eps}{\eps}$ for $\eps>0$ be the first variation of the state variable $u_\eps$. We assume that $\nabla K_\eps \to \nabla K$ and $\eps K_\eps \to 0$, where $K$ solves \eqref{eq:K_limit}. 
        \item[(iii)] For all $\eps >0$ we have $\psi \in \VE(\Dsf_\eps)$ if and only if $\psi\circ T_\eps \in \VE(\Dsf)$.
    \end{itemize}
\end{assumption}
Rearranging and replacing $\psi \circ T_\eps$ by $\psi$ (using Assumption~\ref{ass:main},item (iii)) yields
\begin{equation}\label{eq:Keps}
    \begin{split}
    \int_{\Dsf_\eps} & [  A_1^{(\eps)}(T_\eps(x),u_0 \circ T_\eps + \eps K_\eps, (\nnabla u_0)\circ T_\eps + \nnabla K_\eps)  - A_1^{(\eps)}(T_\eps(x),u_0\circ T_\eps, (\nnabla u_0)\circ T_\eps)] \cdot \eps \psi \\
    +& [  A_2^{(\eps)}(T_\eps(x),u_0 \circ T_\eps + \eps K_\eps, (\nnabla u_0)\circ T_\eps + \nnabla K_\eps)  - A_2^{(\eps)}(T_\eps(x),u_0\circ T_\eps, (\nnabla u_0)\circ T_\eps)] : \nnabla \psi \; dx\\
    =&  \int_{\omega} ((F_1^{\inn}-F_1^{\out})\circ T_\eps ) \cdot ( \eps \psi ) + ((F_2^{\inn} - F_2^{\out})\circ T_\eps ) : \nnabla \psi \; dx  \\
     & -  \int_{\omega} [ A_1^{\inn}-A_1^{\out}](T_\eps(x),u_0\circ T_\eps, (\nnabla u_0)\circ T_\eps) \cdot \eps \psi \; dx\\
     & -  \int_{\omega} [ A_2^{\inn} - A_2^{\out}](T_\eps(x),u_0\circ T_\eps, (\nnabla u_0)\circ T_\eps) :  \nnabla \psi \; dx
\end{split}
\end{equation}
for all $\psi\in \VE(\Dsf_\eps)$. 
With this we can pass to limit $\eps \rightarrow 0$ in \eqref{eq:Keps}  and get 
\begin{align}\label{eq:K_limit} 
    \begin{split}
    \int_{\VR^{d}} & [  A_2^\omega(z,u_0(z), \nnabla u_0(z) + \nnabla K)  - A_2^{\omega}(z,u_0(z), \nnabla u_0(z))] : \nnabla \psi \;dx\\
                   &= \int_{\omega}  (F_2^{\inn} - F_2^{\out})(z) : \nnabla \psi  \;dx  -  \int_{\omega} [ A_2^{\inn} - A_2^{\out}](z,u_0(z), \nnabla u_0(z))  : \nnabla \psi \;dx,
\end{split}
\end{align}
for all $\psi \in \VE(\VR^d):= \lim_{\eps\searrow 0}\VE(\Dsf_\eps)$ , with the definition $A_2^\omega(x, y_1, y_2) = \chi_\omega(x) A_2^{\inn}(x,y_1, y_2) + \chi_{\VR^d \setminus \omega}(x) A_2^{\out}(x,y_1,y_2)$. The limit $\lim_{\eps\searrow 0}\VE(\Dsf_\eps)$ has to be understood formally and in practice often is a Beppo-Levi space (see the next remark).

\begin{remark}\label{rem:beppo_levi}
    The function space $\VE(\VR^d)$ in which \eqref{eq:K_limit} admits a solution can often be chosen as a Beppo-Levi space as shown in \cite{a_ST_2020a}, see also \cite{a_GAST_2020a}. For the numerical approximation in the next section we use precisely that the functions $K_\eps$ converge to $K$ and approximate $K$ by solving an equation on a blown up domain similarly to $\Dsf_\eps$.
\end{remark}

\subsection{Lagrangian theorem for first topological derivatives}
In this section we discuss a method proposed by M.C. Delfour in \cite[Thm.3.3]{c_DE_2018b}. The definite advantage is that it uses the unperturbed adjoint 
equation and only requires the asymptotic analysis of the state equation, but it seems to come with the shortcoming that certain cost functions cannot be treated rigorously; see \cite{a_GAST_2020a} and also \cite{a_BAST_2021a}.

Let $\VE(\Dsf)$ be a Banach space of functions on $\Dsf$. For all parameter $\eps\ge 0$ small consider 
a function $u_\eps\in \VE(\Dsf)$ solving the variational problem of the form
\begin{equation}\label{eq:abstract_state}
    a_\eps[u_\eps](\varphi) = f_\eps(\varphi) \quad \text{ for all } \varphi\in \VE(\Dsf),
\end{equation}
where $a_\eps: \VE(\Dsf) \rightarrow \VE(\Dsf)'$ is a nonlinear operator and $f_\eps$ is a linear form on $\VE(\Dsf)$, respectively. Throughout we assume that this abstract state equation admits a unique solution and 
that $u_\eps-u_0\in \VE(\Dsf)$ for all $\eps$.  Consider now a cost function 
\begin{equation}
    j(\eps) = J_\eps(u_\eps) \in \VR,
\end{equation}
where for all $\eps\ge0$ the functional $J_\eps:\VE(\Dsf)\to \VR$ is differentiable at $u_0$. In the following sections we review methods how to obtain an asymptotic expansion of $j(\eps)$ at $\eps =0$. For this purpose we introduce the Lagrangian function 
\[ 
    \FL(\eps,u,v)= J_\eps(u) + a_\eps[u](v) - f_\eps(v), \quad u \in  \VE(\Dsf), \; v\in \VE(\Dsf).
\]
\begin{proposition}[\cite{c_DE_2018b}]\label{prop:main_df}
    Let $\ell:[0,\tau]\to\VR$ be a function with $\ell(\eps)>0$ for $\eps>0$ and $\displaystyle\lim_{\eps\searrow 0}\ell(\eps)=0$. Furthermore, assume that the limits
\begin{align}
    \FR_1(u_0, p_0) :=& \underset{\eps \searrow 0}{\lim} \; \frac{1}{\ell(\eps)} \int_0^1   \left[\partial_u\FL(\eps,su_\eps+(1-s)u_0, p_0) -  \partial_u  \FL(\eps, u_0, p_0)(u_\eps - u_0)\right]\; ds,\label{eq:deriv_df_11}\\
        \FR_2(u_0, p_0) :=& \underset{\eps \searrow 0}{\lim}\; \frac{1}{\ell(\eps)}(\partial_u \FL(\eps, u_0,p_0) - \partial_u \FL(0,u_0,p_0))(u_\eps - u_0),\label{eq:deriv_df_12}\\
        \partial_{\ell}\FL(0, u_0, p_0) :=&  \lim_{\eps\searrow 0} \;\frac{1}{\ell(\eps)}(\FL(\eps,u_0, p_0) - \FL(0,u_0, p_0)),\label{eq:deriv_df_13}
    \end{align}
    exist. Then the following expansion holds:
    \begin{align} \label{eq_dJ1_direct}
        j(\eps) = j(0) +  \ell(\eps)((\FR_1(u_0, p_0) + \FR_2(u_0, p_0) + \partial_{\ell}\FL(0, u_0, p_0)) + o(\ell(\eps)).
    \end{align}
\end{proposition}

\begin{remark}
Note that by the fundamental theorem of calculus $\FR_1$ can be equivalently written as 
\[
    \FR_1(u_0, p_0) = \underset{\eps \searrow 0}{\lim} \; \frac{1}{\ell(\eps)}  \left[\FL(\eps,u_\eps, p_0) - \FL(\eps, u_0,p_0) -  \partial_u  \FL(\eps, u_0, p_0)(u_\eps - u_0)\right].
\]
However, the form stated in the proposition is better suited for the formal computation of the topological derivatives presented later on.
\end{remark}


\subsection{Formal computation of topological derivative}
We are now applying formally Proposition~\ref{prop:main_df} to compute the topological derivative of problem \eqref{eq_abstrProblem}. For this we let $\FL(\eps,\varphi,\psi):= G(\eps, \varphi, \nnabla \varphi, \psi, \nnabla\psi)$. For sake of better readability, from now on we drop the dependence of $A_i^{(\eps)}$, $i=1,2$ and $j^{(\eps)}$ on the space variable $x$. In what follows we compute the terms $\FR_1(u_0, p_0)$, $\FR_2(u_0, p_0)$ and $\partial_{\ell}\FL(0, u_0, p_0)$  of Proposition~\ref{prop:main_df} for the Lagrangian $\FL$ separately as follows. For this purpose we first introduce the following abbreviations:
\begin{equation*}
\begin{aligned}[c]
    a_1(\eps, \varphi,\psi)  & := \langle A_1^{(\eps)}(\varphi, \nnabla\varphi), \psi \rangle, \\
    f_1(\eps,\psi)  & := \langle F_1^{(\eps)}, \psi \rangle, \\
    J_1(\eps,\psi)  & := J^{\text{vol}}(\Omega_\eps, \psi, \nnabla \psi),
\end{aligned}
\begin{aligned}[c]
    a_2(\eps, \varphi,\psi)  & := \langle A_2^{(\eps)}(\varphi, \nnabla\varphi), \nnabla \psi \rangle,\\
    f_2(\eps,\psi) &  := \langle F_2^{(\eps)}, \nnabla \psi \rangle, \\
    J_2(\psi) &  := J^{\text{bnd}}(\psi, \nnabla \psi).
\end{aligned}
\end{equation*}
Then we split for instance $\FR_1(u_0, p_0)$ into the parts coming from $a_1,a_2,f_1,f_2,J_1,J_2$. The term coming from $a_1$ contributing to $\FR_1(u_0,p_0)$ is for instance given by 
\begin{equation}
    \FR_1^{A_1}(u_0, p_0) :=   \lim_{\eps\searrow 0} \frac{1}{|\omega_\eps|}\int_0^1 \left[\partial_{ \varphi } a_1(\eps,u_0 + s (u_\eps - u_0),p_0) - \partial_{\varphi }a_1(\eps, u_0, p_0)\right](u_\eps-u_0)\; ds,
\end{equation}
where $\partial_{\varphi}a_1$ denotes the derivative with respect to the second argument.
Similarly we denote by $\FR_1^{A_2}(u_0, p_0)$ and $\FR_1^{J_1}(u_0,p_0)$ the contributions of $a_2$ and $J_1$ to the term $\FR_1(u_0,p_0)$, respectively. We proceed in the same fashion for each term. The detailed computations are outlined in the following subsections.
Additionally to Assumption~\ref{ass:main} we make the following assumption.
\begin{assumption}\label{ass:main2}
       We assume that $p_0$ is continuously differentiable at $z$.
\end{assumption}

\subsubsection{Terms coming from $A_1^{(\eps)}$}
For $A_1^{(\eps)} : \VR^m \times \VR^{m \times d} \rightarrow \VR^m$ let $\partial_{y_1} A_1^{(\eps)} : \VR^m \times \VR^{m \times d}\rightarrow \mathcal L(\VR^m, \VR^m)$ denote the derivative with respect to the first argument and $\partial_{y_2} A_1^{(\eps)} : \VR^m \times \VR^{m \times d} \rightarrow \mathcal L( \VR^{m \times d}, \VR^m)$ the derivative with respect to the second argument.
We compute
\begin{align*}
    & \FR_{1, \eps}^{A_1} :=   \int_0^1 \left[\partial_{\varphi} a_1(\eps,u_0 + s (u_\eps - u_0),p_0) - \partial_{\varphi}a_1(\eps, u_0, p_0)\right](u_\eps-u_0)\; ds   \\
    =& \int_0^1 \int_\Dsf [ \partial_{y_1} A_1^{(\eps)}(u_0 + s (u_\eps - u_0), \nnabla u_0 + s \nnabla (u_\eps - u_0) ) -  \partial_{y_1} A_1^{(\eps)}(u_0 , \nnabla u_0  )  ](u_\eps - u_0) \cdot p_0 \; dx \,ds\\
            & + \int_0^1 \int_\Dsf[\partial_{y_2} A_1^{(\eps)}(u_0 + s (u_\eps - u_0), \nnabla u_0 + s \nnabla (u_\eps - u_0) ) \\
            & \hspace{5cm}- \partial_{y_2} A_1^{(\eps)}(u_0 , \nnabla u_0  ) ](\nnabla (u_\eps - u_0) ) \cdot p_0  \; dx \,ds.
\end{align*}
Thus changing variables and using $(\nnabla \varphi)\circ T_\eps = \frac{1}{\eps} \nnabla(\varphi \circ T_\eps)$ and $K_\eps := \frac{(u_\eps - u_0)\circ T_\eps}{\eps}$ yields
\begin{align*}
    \FR_{1,\eps}^{A_1} = &  \eps^d\int_0^1 \int_{\Dsf_\eps} [ \partial_{y_1} A_1^\omega(u_0 \circ T_\eps + s \eps K_\eps, (\nnabla u_0)\circ T_\eps + s \nnabla K_\eps )   \\
                       & \hspace{5cm}-  \partial_{y_1} A_1^\omega(u_0\circ T_\eps , (\nnabla u_0)\circ T_\eps  )  ](\eps K_\eps) \cdot (p_0\circ T_\eps) \; dx \,ds\\
                       &+\eps^d \int_0^1 \int_{\Dsf_\eps}[\partial_{y_2} A_1^\omega(u_0\circ T_\eps + s \eps K_\eps, (\nnabla u_0)\circ T_\eps + s \nnabla K_\eps )  \\
                       & \hspace{5cm}- \partial_{y_2} A_1^\omega(u_0\circ T_\eps , (\nnabla u_0)\circ T_\eps  ) ](\nnabla K_\eps )  \cdot (p_0\circ T_\eps)  \; dx \,ds.
\end{align*}
Therefore setting $\FR_1^{A_1}(u_0, p_0)  = \underset{\eps \searrow 0}{\lim} \frac{1}{|\omega_\eps|} \FR_{1,\eps}^{A_1}$, we obtain
\begin{align*}
    \FR_1^{A_1}(u_0, p_0)  =    \frac{1}{|\omega|} \int_{\VR^d}[A_1^\omega(u_0(z), \nnabla u_0(z) + \nnabla K )   & - A_1^\omega(u_0(z), \nnabla u_0(z) ) \\
    & - \partial_{y_2} A_1^\omega(u_0(z) , \nnabla u_0(z)  ) (\nnabla K )] \cdot p_0(z)  \; dx.
\end{align*}
Next, we compute
\begin{align*}
    \FR_{2,\eps}^{A_1}  := &    \left[\partial_{ \varphi } a_1(\eps,u_0,p_0) - \partial_{ \varphi }a_1( 0 , u_0,p_0)\right](u_\eps-u_0) \\
                      =&   \int_\Dsf [ \partial_{y_1} A_1^{(\eps)}(u_0 , \nnabla u_0  ) -  \partial_{y_1} A_1^{(0)}(u_0 , \nnabla u_0  )  ](u_\eps - u_0) \cdot p_0 \; dx \\
    &+\int_\Dsf[\partial_{y_2} A_1^{(\eps)}(u_0 , \nnabla u_0  ) - \partial_{y_2} A_1^{(0)}(u_0 , \nnabla u_0  ) ](\nnabla (u_\eps - u_0) ) \cdot p_0  \; dx.
\end{align*}
Thus changing variables and using $(\nnabla \varphi)\circ T_\eps = \frac{1}{\eps} \nnabla(\varphi \circ T_\eps)$ and $K_\eps := \frac{(u_\eps - u_0)\circ T_\eps}{\eps}$ yields
\begin{align*}
    \FR_{2,\eps}^{A_1}  = &  \eps^d\int_\omega [ \partial_{y_1} A_1^{\inn}(u_0\circ T_\eps , (\nnabla u_0)\circ T_\eps  ) -  \partial_{y_1} A_1^{\out}(u_0\circ T_\eps , (\nnabla u_0)\circ T_\eps  )  ](\eps K_\eps)  \cdot (p_0\circ T_\eps) \; dx \\
                     &+ \eps^d\int_\omega[\partial_{y_2} A_1^{\inn}(u_0\circ T_\eps , (\nnabla u_0)\circ T_\eps  ) - \partial_{y_2} A_1^{\out}(u_0\circ T_\eps , (\nnabla u_0)\circ T_\eps  ) ](\nnabla K_\eps) \cdot (p_0\circ T_\eps)  \; dx.
\end{align*}
Therefore setting $\FR_2^{A_1}(u_0, p_0)  = \underset{\eps \searrow 0}{\lim} \frac{1}{|\omega_\eps|} \FR_{1,\eps}^{A_1}$, we obtain
\begin{align*}
    &\FR_2^{A_1}(u_0, p_0) = \frac{1}{|\omega|} \int_\omega[\partial_{y_2} A_1^{\inn}(u_0(z) , \nnabla u_0(z)  ) - \partial_{y_2} A_1^{\out}(u_0(z) , \nnabla u_0(z)  ) ](\nnabla K) \cdot p_0(z)  \; dx.
\end{align*}
Finally we compute the term contributing to $\partial_\ell \mathfrak{L}(0, u_0, p_0)$ coming from $A_1^{(\eps)}$ by
\begin{align*} 
    \partial_\ell \mathfrak{L}^{A_1}(0, u_0, p_0) =  [A_1^{\inn}(u_0(z), \nnabla u_0(z)) - A_1^{\out}(u_0(z), \nnabla u_0(z))] \cdot p_0(z).
\end{align*}


\subsubsection{Terms coming from $A_2^{(\eps)}$}
For $A_2^{(\eps)} : \VR^m \times \VR^{m \times d} \rightarrow \VR^{m\times d}$ let $\partial_{y_1} A_2^{(\eps)} : \VR^m \times \VR^{m \times d}\rightarrow \mathcal L(\VR^m, \VR^{m\times d})$ denote the derivative with respect to the first argument and $\partial_{y_2} A_2^{(\eps)} : \VR^m \times \VR^{m \times d} \rightarrow \mathcal L( \VR^{m \times d}, \VR^{m\times d})$ the derivative with respect to the second argument. We have
\begin{align*}
    \FR_{1,\eps}^{A_2}   =& \int_0^1  \left(\partial_{\varphi}a_2(\eps,u_0 + s(u_\eps-u_0),p_0) - \partial_{\varphi}a_2(\eps, u_0, p_0)\right)(u_\eps-u_0)\;ds  \\
                       =& \int_0^1 \int_\Dsf [ \partial_{y_1} A_2^{(\eps)}(u_0 + s (u_\eps - u_0), \nnabla u_0 + s \nnabla (u_\eps - u_0) ) \\
                     & \hspace{6cm} -  \partial_{y_1} A_2^{(\eps)}(u_0 , \nnabla u_0  )  ](u_\eps - u_0) : \nnabla p_0 \; dx \,ds\\
    &+\int_0^1 \int_\Dsf[\partial_{y_2} A_2^{(\eps)}(u_0 + s (u_\eps - u_0), \nnabla u_0 + s \nnabla (u_\eps - u_0) ) \\
    & \hspace{6cm} - \partial_{y_2} A_2^{(\eps)}(u_0 , \nnabla u_0  ) ](\nnabla (u_\eps - u_0) ) : \nnabla p_0  \; dx \,ds.
\end{align*}
Letting $\FR_1^{A_2}(u_0, p_0)  := \underset{\eps \searrow 0}{\lim} \frac{1}{|\omega_\eps|} \FR_{1,\eps}^{A_2}$ and following the same steps as in the computation of $\FR_1^{A_1}(u_0,p_0)$ above leads to 
\begin{equation*}
    \begin{split}
    \FR_1^{A_2}(u_0, p_0) =  \frac{1}{|\omega|} \int_{\VR^d}[A_2^\omega(u_0(z), \nnabla u_0(z) + \nnabla K ) & - A_2^\omega(u_0(z), \nnabla u_0(z) )  \\
& - \partial_{y_2} A_2^\omega(u_0(z) , \nnabla u_0(z)  ) (\nnabla K )] : \nnabla p_0(z)  \; dx.
\end{split}
\end{equation*}
We further have 
\begin{align*}
    \FR_{2,\eps}^{A_2}   :=& \left[\partial_{\varphi} a_2(\eps,u_0) - \partial_{\varphi}a_2(  0 , u_0)\right](u_\eps-u_0) \\
                       = & \int_\Dsf [ \partial_{y_1} A_2^{(\eps)}(u_0 , \nnabla u_0  ) -  \partial_{y_1} A_2^{(0)}(u_0 , \nnabla u_0  )  ](u_\eps - u_0) :  \nnabla p_0 \; dx \\
    &+ \int_\Dsf[\partial_{y_2} A_2^{(\eps)}(u_0 , \nnabla u_0  ) - \partial_{y_2} A_2^{(0)}(u_0 , \nnabla u_0  ) ](\nnabla (u_\eps - u_0) ) : \nnabla p_0  \; dx.
\end{align*}
Letting $\FR_2^{A_2}(u_0, p_0)  := \underset{\eps \searrow 0}{\lim} \frac{1}{|\omega_\eps|} \FR_{2,\eps}^{A_2}$ and following the same steps as in the computation of $\FR_2^{A_1}(u_0,p_0)$ leads to 
\begin{align}
    \FR_2^{A_2}(u_0, p_0) & = \frac{1}{|\omega|} \int_\omega[\partial_{y_2} A_2^{\inn}(u_0(z) , \nnabla u_0(z)  ) - \partial_{y_2} A_2^{\out}(u_0(z) , \nnabla u_0(z)  ) ](\nnabla K) :\nnabla p_0 (z)  \; dx.
\end{align}
Finally, the part of $\partial_\ell \FL(0, u_0, p_0)$ coming from $A_2^{(\eps)}$ reads
\begin{align*} 
    \partial_\ell \FL^{A_2}(0, u_0, p_0) = [A_2^{\inn}(u_0(z), \nnabla u_0(z)) - A_2^{\out}(u_0(z), \nnabla u_0(z))] : \nnabla p_0(z).
\end{align*}


\subsubsection{Terms coming from right hand side}
Since the right hand side does not depend on the solution and therefore $\partial_u F = 0$, there is no contribution to $\FR_1$ and $\FR_2$. It only remains
\begin{align} 
    \partial_\ell \FL^F(0, u_0, p_0) =  - (F_1^{\inn} - F_1^{\out}) \cdot p_0(z) - (F_2^{\inn} - F_2^{\out}) : \nnabla p_0(z).
\end{align}
\subsubsection{Terms coming from $J^{\text{vol}}$}
Let $j_{y_1}^{(\eps)}: \VR^m \times \VR^{m \times d} \rightarrow \mathcal L(\VR^m, \VR)$ denote the derivative of $j^{(\eps)}$ with respect to the first argument and $j_{y_2}^{(\eps)} :  \VR^m \times \VR^{m\times d} \rightarrow \mathcal L(\VR^{m\times d}, \VR)$ the derivative with respect to the second argument such that
\begin{align}
    \partial_u J^{\text{vol}}(\Omega_\eps, u, \nnabla u)(\hat u) = \int_\Dsf j_{y_1}^{(\eps)}( u, \nnabla u) (\hat u) + j_{y_2}^{(\eps)}( u, \nnabla u) (\nnabla \hat u) \; dx.
\end{align}
We have
\begin{align*}
    \FR_{1,\eps}^{J^{\text{vol}}} :=&  \int_0^1(\partial_u J_1(\eps, u_0 + s(u_\eps-u_0)) - \partial_u J_1(\eps, u_0))(u_\eps-u_0)\;ds \\
    = &  \int_0^1 [\partial_u J^{\text{vol}}(\Omega_\eps, u_0 + s (u_\eps - u_0), \nnabla u_0 + s \nnabla(u_\eps - u_0) )  \\
      & \hspace{5cm} - \partial_u J^{\text{vol}}(\Omega_\eps, u_0, \nnabla u_0)](u_\eps - u_0) \; ds \\
    = & \int_0^1 \Bigg\lbrace \int_\Dsf [  j^{(\eps)}_{y_1}(u_0 + s(u_\eps - u_0), \nnabla u_0 + s \nnabla(u_\eps - u_0) ) -  j^{(\eps)}_{y_1}( u_0, \nnabla u_0)](u_\eps - u_0) \; dx  \\
     &+ \int_\Dsf [  j^{(\eps)}_{y_2}(u_0 + s(u_\eps - u_0), \nnabla u_0 + s \nnabla(u_\eps - u_0) ) -  j^{(\eps)}_{y_2}(u_0, \nnabla u_0)] (\nnabla(u_\eps - u_0) ) \; dx  \Bigg\rbrace  \; ds.
\end{align*}
Therefore changing variables and using $(\nnabla \varphi)\circ T_\eps = \frac{1}{\eps} \nnabla(\varphi \circ T_\eps)$ and $K_\eps := \frac{(u_\eps - u_0)\circ T_\eps}{\eps}$ yields
\begin{align*}
    \FR_{1,\eps}^{J^{\text{vol}}}  =&  \eps^d\int_0^1 \Bigg\lbrace  \int_{\Dsf_\eps} [  j^{(\eps)}_{y_1}(u_0\circ T_\eps + s \eps K_\eps, (\nnabla u_0)\circ T_\eps + s \nnabla K_\eps )  \\
                                    & \hspace{5cm}-  j^{(\eps)}_{y_1}( u_0\circ T_\eps, (\nnabla u_0)\circ T_\eps)] (\eps K_\eps) \; dx  \\
                 &+ \eps^d\int_{\Dsf_\eps} [  j^{(\eps)}_{y_2}(u_0\circ T_\eps + s \eps K_\eps, (\nnabla u_0)\circ T_\eps + s  \nnabla K_\eps )  \\
                 & \hspace{5cm}-  j^{(\eps)}_{y_2}(u_0\circ T_\eps, (\nnabla u_0)\circ T_\eps)]  (\nnabla K_\eps) \; dx  \Bigg\rbrace  \; ds.
\end{align*}
Letting $\FR_1^{J^{\text{vol}}}(u_0,p_0) := \lim_{\eps\searrow 0} \frac{1}{|\omega_\eps|} \FR_{1,\eps}^{J^{\text{vol}}}$ and using $\eps K_\eps \to 0$ and $\nnabla K_\eps \to \nnabla K$, we obtain 
\begin{align*}
    \FR_1^{J^{\text{vol}}}(u_0,p_0) & = \frac{1}{|\omega|}  \int_{\VR^d} [  j^\omega(u_0(z), \nnabla u_0(z) +  \nnabla K ) -j^\omega(u_0(z), \nnabla u_0(z)  )  \\
    & \hspace{5cm} -  j^\omega_{y_2}(u_0(z), \nnabla u_0(z) )  (\nnabla K) ]\; dx.  
\end{align*}
Here, $j^\omega(y_1, y_2) = \chi_\omega \, j^{\text{in}}(y_1, y_2) + \chi_{\VR^d \setminus \omega} \,j^{\text{out}}(y_1, y_2)$.
\newline
We proceed with the contribution of $J^{\text{vol}}$ to the term $\FR_2(u_0,p_0)$. For this we compute
\begin{align*}
    \FR_{2,\eps}^{J^{\text{vol}}}   :=&  [\partial_u J_1(\eps, u_0)-\partial_u J_1(0, u_0) ](u_\eps - u_0) \\
                                     =&    [ \partial_u J^{\text{vol}}(\eps, u_0, \nnabla u_0) -  \partial_u J^{\text{vol}}(0, u_0, \nnabla u_0)](u_\eps - u_0) \\
                                     =&   \int_\Dsf ( j_{y_1}^{(\eps)}( u_0, \nnabla u_0)   -   j_{y_1}^{(0)}( u_0, \nnabla u_0)) (u_\eps - u_0)  \\
                     & +\int_\Dsf ( j_{y_2}^{(\eps)}( u_0, \nnabla u_0)   -   j_{y_2}^{(0)}( u_0, \nnabla u_0)) (\nnabla (u_\eps - u_0)) .
\end{align*}
Therefore setting $\FR_2^{J^{\text{vol}}}(u_0,p_0):= \lim_{\eps\searrow 0} \frac{1}{|\omega_\eps|} \FR_{2,\eps}^{J^{\text{vol}}}$, we obtain
\begin{align*}
    \FR_2^{J^{\text{vol}}}(u_0,p_0) =&  \underset{\eps \searrow 0}{\lim} \frac{1}{|\omega_\eps|} [ \partial_u J^{\text{vol}}(\eps, u_0, \nnabla u_0) -  \partial_u J^{\text{vol}}(0, u_0, \nnabla u_0)](u_\eps - u_0) \\
    =&  \underset{\eps \searrow 0}{\lim} \frac{1}{|\omega_\eps|} \Bigg[ \int_\Dsf ( j_{y_1}^{(\eps)}( u_0, \nnabla u_0)   -   j_{y_1}^{(0)}( u_0, \nnabla u_0)) (u_\eps - u_0) \\
     &\hspace{2cm} + ( j_{y_2}^{(\eps)}( u_0, \nnabla u_0) -    
 j_{y_2}^{(0)}(u_0, \nnabla u_0) ) (\nnabla (u_\eps - u_0)) \; dx \Bigg] \\
        =&\underset{\eps \searrow 0}{\lim} \frac{1}{|\omega|} \Bigg[ \int_{\omega}[ j_{y_1}^{(\eps)}(u_0\circ T_\eps, (\nnabla u_0)\circ T_\eps) -  j_{y_1}^{(0)}(u_0\circ T_\eps, (\nnabla u_0)\circ T_\eps)](\eps K_\eps) \\
         &\hspace{1cm} + [ j_{y_2}^{(\eps)}(u_0 \circ T_\eps, (\nnabla u_0)\circ T_\eps) -  j_{y_2}^{(0)}(u_0 \circ T_\eps, (\nnabla u_0)\circ T_\eps)](\nnabla K_\eps ) \; dx \Bigg] \\
            =&  \frac{1}{|\omega|}  \int_{\omega} [ j_{y_2}^{\inn}( u_0 (z), \nnabla u_0(z)) - j_{y_2}^{\out}(u_0 (z), \nnabla u_0(z))] (\nnabla K ) \; dx .
\end{align*}
Finally, the term of $J^{\text{vol}}$ contributing to $\partial_\ell \FL$ reads
\begin{align}
    \partial_\ell \FL^{J^{\text{vol}}}(0, u_0, p_0) =  j^{\inn}(u_0(z), \nnabla u_0(z)) - j^{\out}(u_0(z), \nnabla u_0(z)) .
\end{align}

\subsubsection{Terms coming from $J^{\text{bnd}}$}
There is no contribution from the term $J^{\text{bnd}}$, since the rescaled term $\partial \Dsf_\eps$ tends to "infinity" and is 
therefore not present.

\subsubsection{Summary} \label{sec_summary}

In summary we showed that 
\begin{align*}
    \FR_1(u_0,p_0) & =  \FR_1^{A_1}(u_0,p_0) + \FR_1^{A_2}(u_0,p_0)  + \FR_1^{J^{\text{vol}}}(u_0,p_0),\\
    \FR_2(u_0,p_0) & =  \FR_2^{A_1}(u_0,p_0) + \FR_2^{A_2}(u_0,p_0)  + \FR_2^{J^{\text{vol}}}(u_0,p_0), \\
    \partial_\ell \FL(u_0,p_0) & =  \partial_\ell \FL^{A_1}(u_0,p_0) + \partial_\ell \FL^{A_2}(u_0,p_0) + \partial_\ell \FL^{F}(u_0,p_0) + \partial_\ell \FL^{J^{\text{vol}}}(u_0,p_0).
\end{align*}
Therefore we proved the following theorem under the assumption that formally $\nabla K_\eps \to \nabla K$ in $\VE(\VR^d)$ and $\eps K_\eps \rightarrow 0$ and that $u_0$ is sufficiently smooth at $z$. 

\begin{theorem} \label{thm_formulaTD}
The first order topological expansion of $\mathcal J$ at $\Omega)$ in the points $z$ reads
\begin{align} \label{eq_summary_dJ}
    d \mathcal J(\Omega, \omega)(z) = \FR_1(u_0, p_0) + \FR_2(u_0,p_0) + \partial_\ell \FL(0, u_0, p_0),
\end{align}
where
\begin{align*}
    \FR_1(u_0, p_0) =& \frac{1}{|\omega|} \int_{\VR^d}[A_1^\omega(u_0(z), \nnabla u_0(z) + \nnabla K ) - A_1^\omega(u_0(z), \nnabla u_0(z) )  \\
                   & \hspace{5cm} - \partial_{y_2} A_1^\omega(u_0(z) , \nnabla u_0(z)  ) (\nnabla K )] \cdot p_0(z)  \; dx \\
                   &+\frac{1}{|\omega|} \int_{\VR^d}[A_2^\omega(u_0(z), \nnabla u_0(z) + \nnabla K ) - A_2^\omega(u_0(z), \nnabla u_0(z) )  \\
                   & \hspace{5cm}  - \partial_{y_2} A_2^\omega(u_0(z) , \nnabla u_0(z)  ) (\nnabla K )] : \nnabla p_0(z)  \; dx \\
                   &+\frac{1}{|\omega|}  \int_{\VR^d} [  j^\omega(u_0(z), \nnabla u_0(z) +  \nnabla K(x) ) -j^\omega(u_0(z), \nnabla u_0(z)  )  \\
                   & \hspace{6cm} -  j^\omega_{y_2}(u_0(z), \nnabla u_0(z) ) ( \nnabla K) ]\; dx \\     
    \FR_2(u_0, p_0) =&\frac{1}{|\omega|} \int_\omega[\partial_{y_2} A_1^{\inn}(u_0(z) , \nnabla u_0(z)  ) - \partial_{y_2} A_1^{\out}(u_0(z) , \nnabla u_0(z)  ) ](\nnabla K) \cdot p_0(z)  \; dx \\
                   &+\frac{1}{|\omega|} \int_\omega[\partial_{y_2} A_2^{\inn}(u_0(z) , \nnabla u_0(z)  ) - \partial_{y_2} A_2^{\out}(u_0(z) , \nnabla u_0(z)  ) ](\nnabla K) : \nnabla p_0 (z)  \; dx \\
                   &+ \frac{1}{|\omega|}  \int_{\omega} [ j_{y_2}^{\inn}( u_0 (z), \nnabla u_0(z)) - j_{y_2}^{\out}( u_0 (z), \nnabla u_0(z))] (\nnabla K ) \; dx \\
                   \intertext{ and } 
    \partial_\ell \FL(0, u_0, p_0) =& [A_1^{\inn}(u_0(z), \nnabla u_0(z)) - A_1^{\out}(u_0(z), \nnabla u_0(z))] \cdot p_0(z) \\
                                  &+  [A_2^{\inn}(u_0(z), \nnabla u_0(z)) - A_2^{\out}(u_0(z), \nnabla u_0(z))] : \nnabla p_0(z)  \\
                                  &- (F_1^{\inn} - F_1^{\out}) \cdot p_0(z) - (F_2^{\inn} - F_2^{\out}) : \nnabla p_0(z) \\
                                  &+  j^{\inn}( u_0(z), \nnabla u_0(z)) - j^{\out}(u_0(z), \nnabla u_0(z)).
\end{align*}
The function $K$ solves
\begin{align}  \label{eq_K}
    \begin{aligned}
    \int_{\VR^{d}} & [  A_2^\omega(z,u_0(z), \nnabla u_0(z) + \nnabla K)  - A_2^{\omega}(z,u_0(z), \nnabla u_0(z))] : \nnabla \psi \; dx\\
                   &= \int_{\omega}  (F_2^{\inn} - F_2^{\out})(z) : \nnabla \psi  \; dx  -  \int_{\omega} [ A_2^{\inn} - A_2^{\out}](z,u_0(z), \nnabla u_0(z))  : \nnabla \psi\; dx
    \end{aligned}
\end{align}
for all $\psi \in \VE(\VR^d)$.
\end{theorem}

\section{Numerical experiments} \label{sec_numExp}
In this section, we consider several linear and nonlinear problems in two and three space dimensions for which we numerically compute the topological derivative by means of the procedure outlined in the previous section, see the generic formulas of Theorem \ref{thm_formulaTD}. We solve all underlying PDEs by means of piecewise linear and globally continuous finite elements on triangular grids using the finite element software package \texttt{NGSolve} \cite{Schoeberl2014}. We compute an approximation of the corrector function $K$, which is defined as the solution of a PDE on the unbounded domain $\VR^d$, by solving \eqref{eq_K} on a large, but bounded domain $B_R:= B(\mathbf 0, R)$ with homogeneous Dirichlet boundary conditions on $\partial B_R$. This procedure is motivated by the fact that the solution $K$ often can be shown to exhibit a decay behaviour as $|x| \rightarrow \infty$. We will show that the numerically computed topological derivative matches the analytical formula well for the topological derivative in cases where this formula is known. Moreover, we perform numerical tests to verify the topological asymptotic expansion
\begin{align} \label{eq_topAsympExp}
    \mathcal J(\Omega_\eps) = \mathcal J(\Omega) + \eps^d |\omega| d \mathcal J(\Omega)(z) + \mathcal O(\eps^{d+1}),
\end{align}
or, in other words
\begin{align} \label{eq_deltaJ}
    \delta \mathcal J := |\mathcal J(\Omega_\eps) - \mathcal J(\Omega) + \eps^d |\omega| d \mathcal J(\Omega)(z) | = \mathcal O(\eps^{d+1}).
\end{align}
We consider several different inclusion shapes in two and three space dimensions. In 2D, we consider
\begin{enumerate}
    \item $\omega = \omega_{2D}^{(1)} := B( (0,0), 1)$ the unit disk in $\VR^2$
    \item $\omega = \omega_{2D}^{(2)} := B( (0.5,0.5), 1)$ a shifted unit disk in $\VR^2$
    \item $\omega = \omega_{2D}^{(3)} := Ell( (0,0), \frac{3}{2}, \frac{2}{3})$ an axis-aligned ellipse centered at the origin with axes lengths $\frac32$ and $\frac23$,
    \item $\omega = \omega_{2D}^{(4)} := Ell( (0.5,0.5), \frac{3}{2}, \frac{2}{3})$ a shifted axis-aligned ellipse centered at the point $(0.5,0.5)^\top$ with the same axes lengths $\frac32$ and $\frac23$,
    \item $\omega = \omega_{2D}^{(5)} :=$ an L-shaped domain of area $\pi$,  
\end{enumerate}
and in 3D the shapes
\begin{enumerate}
    \item $\omega = \omega_{3D}^{(1)} := B( (0,0,0), 1)$ the unit ball in $\VR^3$
    \item $\omega = \omega_{3D}^{(2)} := B( (0.5,0.5,0.5), 1)$ a shifted unit ball in $\VR^3$
    \item $\omega = \omega_{3D}^{(3)} := Ell( (0,0,0), \frac{3}{2}, \frac{2}{3}, 1)$ an axis-aligned ellipsoid centered at the origin with axes lengths $\frac32$, $\frac23$ and $1$,
    \item $\omega = \omega_{3D}^{(4)} := Ell( (0.5,0.5,0.5), \frac{3}{2}, \frac{2}{3}, 1)$ a shifted axis-aligned ellipsoid centered at the point $(0.5, 0.5, 0.5)^\top$ with the same axes lengths $\frac32$, $\frac23$ and $1$.
\end{enumerate}
Note that $\omega_{2D}^{(1)}$, $\omega_{2D}^{(3)}$, $\omega_{3D}^{(1)}$ and $\omega_{3D}^{(3)}$ are symmetric with respect to the $x_1$-, the $x_2$- (and the $x_3$-)axes while the other shapes are not.

The finite element software $\texttt{NGSolve}$ allows to define PDEs in weak form in a symbolic way and also supports automated differentiation of expressions, see \cite{GanglSturmNeunteufelSchoeberl2020} for applications of these capabilities in the context of shape derivatives. Our implementation is available from \cite{GanglSturmCodeTDAuto}. It consists of a main file, which is completely independent of the concrete topology optimization problem at hand, and four other files defining the geometry, the PDE, the cost function as well as some algorithmic parameters. The main file implements the topological derivative of a general problem of the form \eqref{eq_abstrProblem}, i.e., it implements the solution of the state and adjoint equation, of the corresponding corrector equation \eqref{eq_K} on the large domain $B_R$ as well as the computation of the terms $\FR_1(u_0, p_0)$, $\FR_2(u_0, p_0)$ and $\partial_\ell \FL(0,u_0, p_0)$ as they are given in Section \ref{sec_summary}.

\subsection{Diffusion, convection, reaction}
Here, we consider the class of topology optimization problems with tracking-type cost functionals and a scalar diffusion-convection-reaction equation as a PDE constraint. Given the computational domain whose boundary is divided into the disjoint Dirichlet and Neumann boundaries, $\partial \Dsf = \Gamma_D \cup \Gamma_N$, the PDE-constrained topology optimization problem is to find $u \in H^1_{\Gamma_D}(\Dsf) := \{v \in H^1(\Dsf): v|_{\Gamma_D}=0 \}$ and $\Omega \in \mathcal A$ for some set of admissible shapes $\mathcal A$ as a solution to
\begin{align} \label{eq_exH1_J}
  \underset{\Omega}{\mbox{min }} J(u, \Omega):= 
  \int_\Dsf \tilde \alpha_\Omega(x) |u-u_d|^2 + \tilde \beta_\Omega(x) |\nabla(u-u_d)|^2 \; \mbox dx +&  \tilde \gamma \int_{\Gamma_N} |u - u_d|^2 \; \mbox ds
  \end{align}
  such that 
  \begin{align}
      \int_\Dsf \beta_\Omega(x,|\nabla u|) \nabla u & \cdot \nabla \psi + (\mathbf{b}_\Omega(x)\cdot \nabla u) \psi + \alpha_\Omega(x, u) \psi \; \mbox dx  \\
                                                                                                      =& \int_\Dsf f_\Omega(x) \psi + \mathbf M_\Omega(x) \cdot \nabla \psi \; \mbox dx + \int_{\Gamma_N} g_N \cdot \psi \; \mbox dS \label{eq_exH1_pde}
\end{align}
for all $\psi \in H^1_{\Gamma_D}(\Dsf)$. Here, $c_\Omega(x) := \chi_{\Omega}(x) c_{1} + \chi_{\Dsf \setminus \Omega}(x) c_{2}$ for some given constants $c_{1}$, $c_{2}$ for $c \in \{\tilde \alpha, \tilde \beta, \mathbf{b}, f, \mathbf M \}$. The functions $\beta_\Omega(x, |\nabla u|)$ and $\alpha_\Omega(x, u)$ will be defined piecewise in the subsequent subsections.

This problem fits into the framework considered in Section \ref{sec_derivTD} with the choices
\begin{align*}
    j^\Omega(x, u, \nabla u) =& \tilde \alpha_\Omega(x) |u-u_d|^2 + \tilde \beta_\Omega(x) |\nabla(u-u_d)|^2, \\
    j^{bnd}(x, u, \nabla u) =&  \tilde \gamma |u-u_d|^2, \\
    A_1^\Omega(x, u, \nabla u) =& (\mathbf b_\Omega(x) \cdot \nabla u) + \alpha_\Omega(x, u), \\
    A_2^\Omega(x, u, \nabla u) =& \beta_\Omega(x, |\nabla u|)\nabla u \\
    F_1^\Omega(x) =& f_\Omega(x), \\
    F_2^\Omega(x) =& \mathbf M_\Omega(x) .
\end{align*}

\subsubsection{Example 1: A linear diffusion-convection-reaction problem in 2D}   \label{sec_2dlin}
We begin with a simple two-dimensional, linear version of problem \eqref{eq_exH1_J}--\eqref{eq_exH1_pde} where we set $u_d = 0$, $\tilde \alpha_1 =1$, $\tilde \alpha_2 =2$, $\tilde \beta_1 = \tilde \beta_2 =0$, $\tilde \gamma = 0$, $\mathbf b_1 = (1,0)^\top$, $\mathbf b_2 = (0,1)^\top$, $f_1 = 1$, $f_2=2$, $\mathbf M_1 = \mathbf M_2 = (0,0)^\top$, $g_N(x_1, x_2) = x_1 x_2$. 
Moreover, we choose the piecewise constant functions
\begin{align*}
    \beta_\eps(x, |\nabla u|) =& \chi_{\Omega_\eps}(x) \beta_1 + \chi_{\Dsf \setminus \Omega_\eps}(x) \beta_2, \\
    \alpha_\eps(x,  u) =& \chi_{\Omega_\eps}(x) \alpha_1 + \chi_{\Dsf \setminus \Omega_\eps}(x) \alpha_2, 
\end{align*}
with the values $\beta_1 = 1$, $\beta_2 = 2$, $\alpha_1 = 1$, $\alpha_2 = 2$, which thus do not depend on $u$ or $\nabla u$.

In this setting, it is well-known that, for a point $z \in \Dsf \setminus \Omega$ and inclusion shape $\omega = B_1(0)$ the unit disk, the corrector term $K$ is linear inside $\omega$ and satisfies
\begin{align}
    \nabla K|_{\omega} = - \frac{\beta_1 - \beta_2}{\beta_1+\beta_2} \nabla u(z),
\end{align}
see e.g. \cite{gangl2020topological}. Thus, it follows that
\begin{align}
  \FR_2^{A_2}(u,p) =& \frac{1}{|\omega|} \int_\omega (\beta_1  - \beta_2)  \nabla K \cdot \nabla p(z) \; dx 
  = -  (\beta_1  - \beta_2)   \frac{\beta_1 - \beta_2}{\beta_1+\beta_2} \nabla u(z) \cdot \nabla p(z) 
\end{align}
which together with the term $\partial_\ell \FL^{A_2}(u, p) = (\beta_1 - \beta_2) \nabla u(z) \cdot \nabla p(z)$ sums up to
\begin{align}
 \FR_2^{A_2}(u,p) +\partial_\ell \FL^{A_2}(u, p)
 =& 2\beta_2 \frac{\beta_1 - \beta_2}{\beta_1+\beta_2}  \nabla u(z) \cdot \nabla p(z).
\end{align}
Similarly, $\partial_\ell \FL^{A_1}(u,p) = (\mathbf b_1 - \mathbf b_2) \cdot \nabla u(z) \, p(z) + (\alpha_1 - \alpha_2)u(z)p(z)$ together with
\begin{align}
    \FR_2^{A_1}(u,p) =& \frac{1}{|\omega|} \int_\omega (\mathbf b_1  - \mathbf b_2) \cdot \nabla K \, p(z) \; dx 
    = -  \frac{\beta_1 - \beta_2}{\beta_1+\beta_2}   (\mathbf b_1  - \mathbf b_2) \cdot \nabla u(z) \, p(z) 
\end{align}
adds up to
\begin{align}
    \FR_2^{A_1}(u,p) + \partial_\ell \FL^{A_1}(u,p) = 2 \frac{\beta_2}{\beta_1+\beta_2} (\mathbf b_1  - \mathbf b_2) \cdot \nabla u(z) \, p(z)  +  (\alpha_1 - \alpha_2)u(z)p(z).
\end{align}
Moreover, we have $\FR_1(u_0, p_0) = 0$ since the PDE constraint is linear and the cost function does not depend on $\nabla u$. Together with the remaining terms $\partial_\ell \FL^{F}$ and $\partial_\ell \FL^{J^{\text{vol}}}$,
the topological derivative reads in closed form,
\begin{align} \label{eq_TD_Ex1_closedForm}
    \begin{aligned}
    d \mathcal J(\Omega, \omega)(z) =& 2 \beta_2 \frac{\beta_1 - \beta_2}{\beta_1 + \beta_2} \nabla u (z) \cdot \nabla p(z) +2  \frac{\beta_2}{\beta_1 + \beta_2} (\mathbf b_1 - \mathbf b_2)\cdot \nabla u(z) \, p(z) + (\alpha_1 - \alpha_2)u(z) p(z) \\
    &- (f_1 - f_2)p(z) + (\tilde \alpha_1 - \tilde \alpha_2)u(z)^2,
    \end{aligned}
\end{align}
see also, e.g., \cite{Amstutz2006} for the case where $\mathbf b_1 = \mathbf b_2 = (0,0)^\top$. A similar problem including a convection term is covered by the analysis in \cite{Amstutz_2m_2014}, however, there the factor $2\beta_2 / (\beta_1+\beta_2)$ (and thus the term $\FR_2^{A_1}(u,p)$) is missing. We remark that our numerical experiments indicate that this factor is important for having the correct topological derivative formula.

\begin{figure}
     \begin{tabular}{ccc}
        \includegraphics[width=.35 \textwidth, trim = 80 0 60 0, clip]{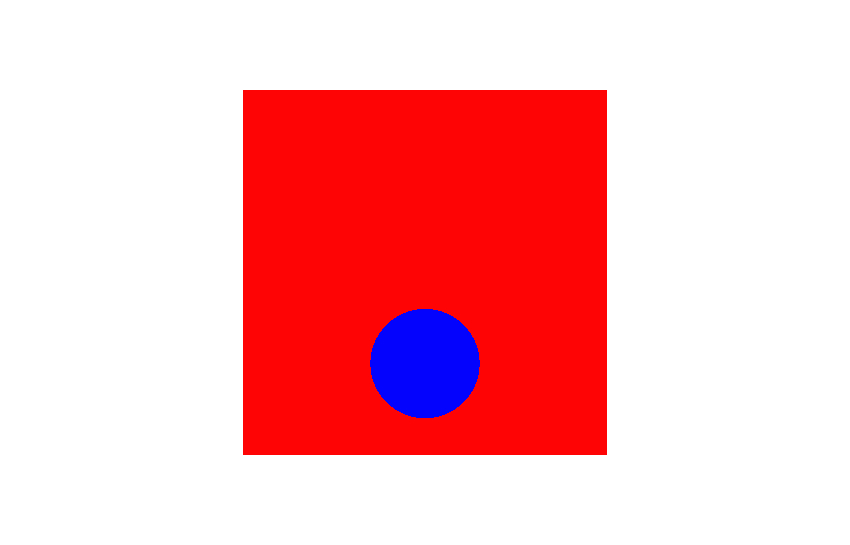}&
        \includegraphics[width=.35 \textwidth, trim = 80 0 60 0, clip]{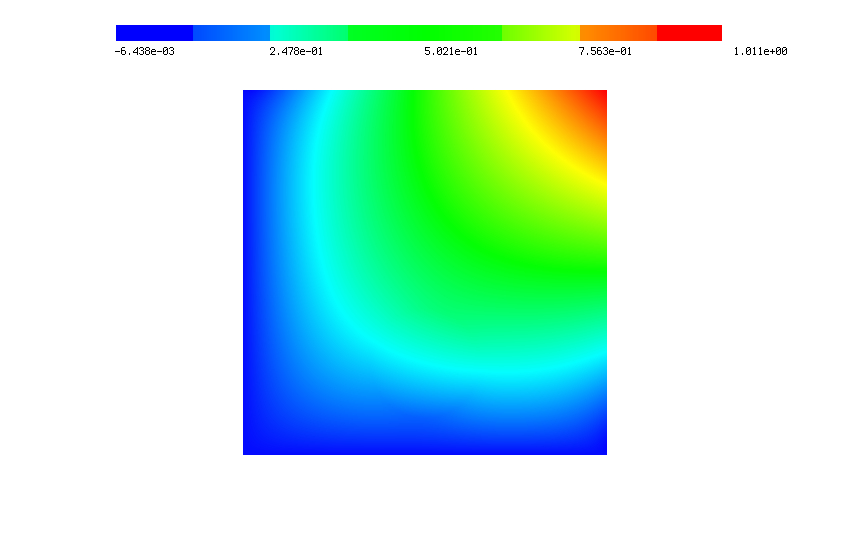} &
        \includegraphics[width=.35 \textwidth, trim = 80 0 60 0, clip]{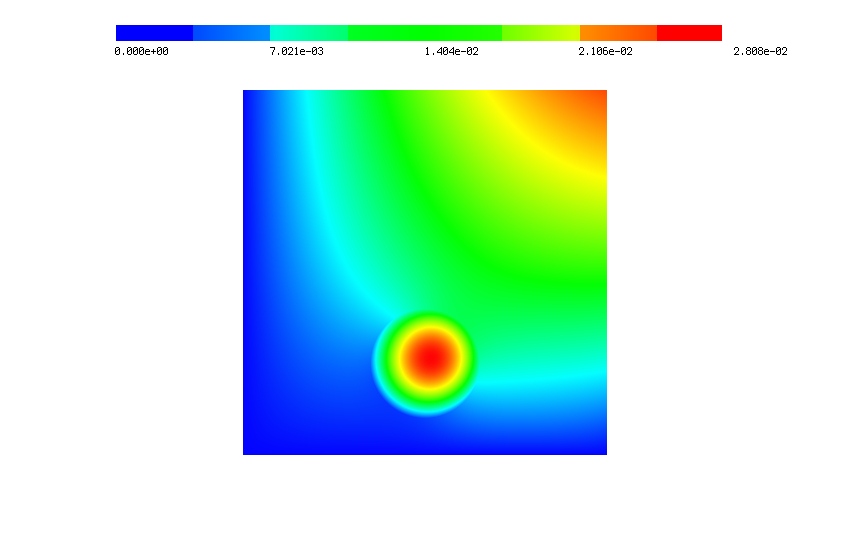} \\
        (a) & (b)&(c)
     \end{tabular}
     \caption{(a) Computational domain $\Dsf$ with subdomain $\Omega$ for Examples 1 and 2. (b) Solution $u$ to unperturbed state equation for linear Example 1. (c) Solution $u$ to unperturbed state equation for nonlinear Example 2.}
     \label{fig_settingEx1}
\end{figure}

We numerically compute the topological derivative by the procedure explained in Section \ref{sec_derivTD} for $\Dsf =(-1,1)^2$ $\Gamma_D = \{(x_1,x_2): x_1 = -1 \mbox{ or } x_2=-1\}$, $\Gamma_N = \partial \Dsf \setminus \Gamma_D$ and $\Omega = B((0,-0.5)^\top, 0.3)$ at the point $z = (0, 0.5)^\top \in \Dsf \setminus \Omega$ for the five inclusion shapes $\omega_{2D}^{(i)}$, $i =1,\dots, 5$ defined above. Figure \ref{fig_settingEx1} depicts the computational domain $\Dsf$ as well as the unperturbed state obtained on a mesh with 528724 vertices. 

\begin{figure}
    \begin{center}\includegraphics[width=.5\textwidth]{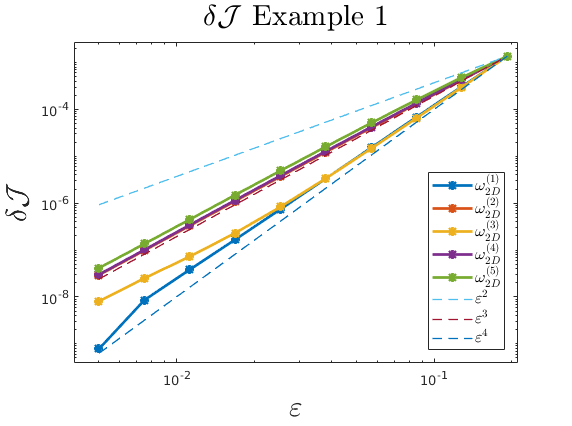}\end{center}
    \caption{Numerical verification of topological asymptotic expansion \eqref{eq_topAsympExp} for Example 1 for five different inclusion shapes.}
    \label{fig_deltaJ_ex1}
\end{figure}

We solved the corresponding problem \eqref{eq_K} for the corrector $K$ on the bounded domain $B(\mathbf 0, R)$ with $R = 1000$ using a finite element mesh consisting of 371950 vertices.
In order to numerically verify the computed values, we performed a Taylor test for the topological asymptotic expansion \eqref{eq_topAsympExp}. We chose a set of decreasing values for the inclusion radius $\eps$, namely 
\begin{equation} \label{eq_defEpsUnderline}
\underline \eps := \{ \eps_0 \, \delta^{9}, \dots, \eps_0 \, \delta^1, \eps_0 \, \delta^0 \} 
\end{equation}
with two constants $\eps_0=0.005$ and $\delta = 1.5$, and solved the perturbed PDE constraint \eqref{eq_pertState} on a mesh which is highly refined around the point $z$ and evaluated the cost function to get perturbed cost function values $\mathcal J(\Omega_\eps)$. For all $\eps \in \underline \eps$ and for each of the five inclusion shapes $\omega_{2D}^{(i)}$, we computed the quantities $\delta \mathcal J$ defined in \eqref{eq_deltaJ}. In Figure \ref{fig_deltaJ_ex1}, it can be observed that, for each inclusion shape, $\delta \mathcal J$ behaves (at least) like $\eps^{3}$. In the case of the two inclusions $\omega_{2D}^{(1)}$, $\omega_{2D}^{(3)}$ which are symmetric with respect to both the $x_1$- and the $x_2$-axes, we can even observe that $\delta \mathcal J = \mathcal O(\eps^{4})$. This is expected since it is known that the second term in the expansion \eqref{eq_topAsympExp} vanishes, $d^2 \mathcal J(\Omega, \omega)(z) = 0$ for this example if symmetric inclusion shapes $\omega$ are considered, see e.g. \cite{a_BAST_2021a}. In Figure \ref{fig_deltaJ_ex1} as well as in the subsequent Taylor test graphs, the data for $\omega_{2D}^{(i)}$, $i>1$, is scaled such that they coincide with the data for $\omega_{2D}^{(1)}$ for the largest considered value of $\eps$. This is done for better comparison of the convergence rates.

\begin{figure}
    \begin{center}\includegraphics[width=.5\textwidth]{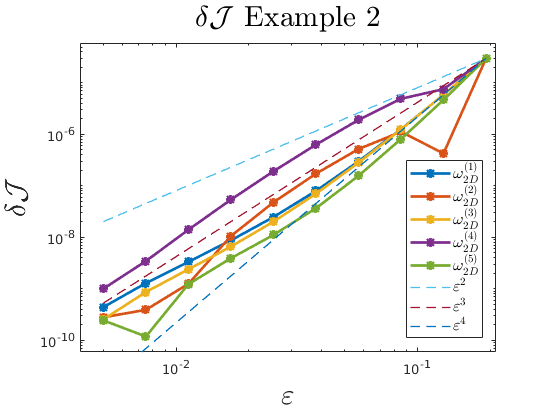}\end{center}
    \caption{Numerical verification of topological asymptotic expansion \eqref{eq_topAsympExp} for Example 2 for five different inclusion shapes.}
    \label{fig_deltaJ_ex2}
\end{figure}

We will further investigate the numerically computed topological derivative for this example in Section \ref{sec_effEval} where we will also discuss the efficient evaluation of the topological derivative in the full computational domain $\Dsf$.

\subsubsection{Example 2: A nonlinear diffusion-convection-reaction problem in 2D}   \label{sec_2d_nonlin}
Next, we consider a more general, quasilinear setting of problem \eqref{eq_exH1_J}--\eqref{eq_exH1_pde} where we set $u_d = 0$, $\tilde \alpha_1 =1$, $\tilde \alpha_2 =2$, $\tilde \beta_1 =1$, $ \tilde \beta_2 =2$, $\tilde \gamma = 1$, $\mathbf b_1 = (1,0)^\top$, $\mathbf b_2 = (0,1)^\top$, $f_1 = 1$, $f_2=2$, $\mathbf M_1 =(1,0)^\top$, $ \mathbf M_2 = (0,1)^\top$, $g_N(x_1, x_2) = x_1 x_2$ as well as the nonlinear functions

\begin{align*}
    \beta_\eps(x, |\nabla u(x)|) =\begin{cases}
         1 & x \in \Omega_\eps, \\
         \hat \beta_2(|\nabla u(x)|) & \mbox{else},
    \end{cases} \qquad
    \alpha_\eps(x, u(x)) =\begin{cases}
         1 & x \in \Omega_\eps, \\
         u(x)^3 & \mbox{else},
    \end{cases}
\end{align*}
where $\hat \beta_2(s) = \nu_0-(\nu_0-200) e^{ - s^6 / 1000}$ with $\nu_0 = 10^7 / (4 \pi)$. This function is sometimes used as a magentic reluctivity function in electromagnetics and satisfies the monotonicity and Lipschitz conditions which ensure existence of a unique solution to the PDE constraint. Topological derivatives for problems involving quasilinear PDE constraints are challenging from both the analytical and the numerical point of view, see e.g. \cite{a_AMBO_2017a, a_AMGA_2019a, a_GAST_2020a, GanglSturm2021Hcurl}.
The numerical experiments depicted in Figure \ref{fig_deltaJ_ex2} seem to exhibit the behavior $\mathcal O(\eps^3)$ for all five inclusion shapes which is in accordance with \eqref{eq_deltaJ} and therefore confirms the numerically computed value for the topological derivative.

\subsubsection{Example 3: A linear diffusion-convection-reaction problem in 3D}   \label{sec_3d_lin}
\begin{figure}
    \begin{center}\includegraphics[width=.5\textwidth]{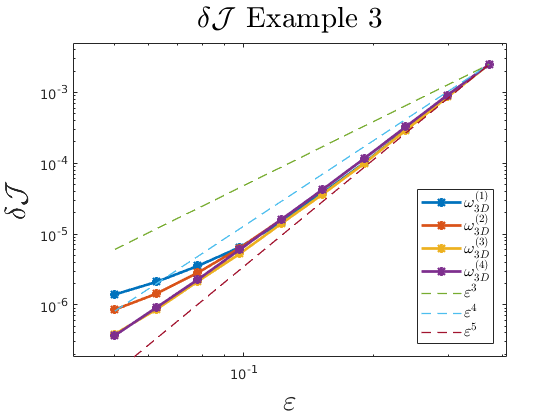}\end{center}
    \caption{Numerical verification of topological asymptotic expansion \eqref{eq_topAsympExp} for Example 3 for four different inclusion shapes.}
    \label{fig_deltaJ_ex3}
\end{figure}

Next, we consider a linear three-dimensional version of \eqref{eq_exH1_J}--\eqref{eq_exH1_pde} on the domain $\Dsf = B(\mathbf 0, 1)$ with the subdomain $\Omega = B((0, -0.7, 0)^\top, 0.2)$. We chose the parameters $u_d = 0$, $\tilde \alpha_1 =1$, $\tilde \alpha_2 =2$, $\tilde \beta_1 = \tilde \beta_2 =0$, $\tilde \gamma = 0$, $\beta_1 = 1$, $\beta_2 = 2$, $\mathbf b_1 = \mathbf b_2 = (0,0)^\top$, $f_1 = 1$, $f_2=2$, $\mathbf M_1 = \mathbf M_2 = (0,0)^\top$, $g_N(x_1, x_2) = x_1 x_2$ and evaluate the topological derivative at the point $z = (0, 0.1, 0)^\top$. For each of the four inclusion shapes $\omega_{3D}^{(i)}$, $i = 1,2,3,4$, we solved the corresponding corrector equation \eqref{eq_K} on a three-dimensional ball $B_R$ of radius $R=1000$ using a tetrahedral mesh with 301116 vertices. Subsequently, we computed the topological derivative according to Section \ref{sec_summary} and computed the quantities $\delta \mathcal J$ \eqref{eq_deltaJ} for the vector of radius values \eqref{eq_defEpsUnderline} with $\eps_0 = 0.05$ and $\delta = 1.25$. Figure \ref{fig_deltaJ_ex3} shows that the quantities $\delta \mathcal J$ decay at least as fast as $\eps^4$ for all four inclusion shapes, which is in accordance with \eqref{eq_topAsympExp}. The deterioration of the rates in Figure \ref{fig_deltaJ_ex3} can be attributed to the discretization error. We used a mesh consisting of 245177 vertices for the computational domain $\Dsf$ which is highly refined around the point $z$.

\begin{figure}
    \begin{center}\includegraphics[width=.5\textwidth]{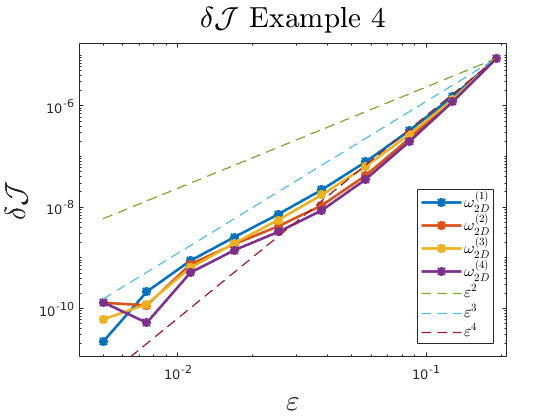}\end{center}
    \caption{Numerical verification of topological asymptotic expansion \eqref{eq_topAsympExp} for Example 4 for four different inclusion shapes.}
    \label{fig_deltaJ_ex4}
\end{figure}

\subsection{Elasticity}
\begin{figure}
    \begin{tabular}{cc}
        \includegraphics[width=.5\textwidth]{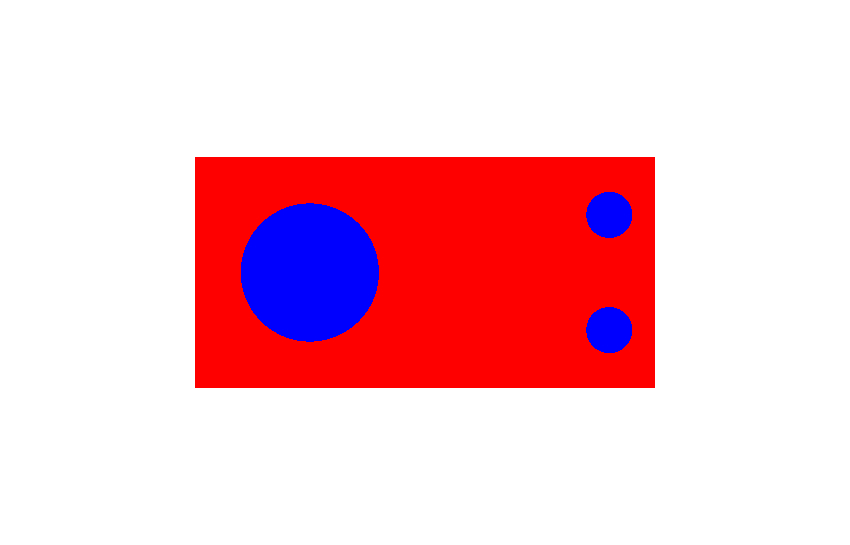} & 
        \includegraphics[width=.5\textwidth]{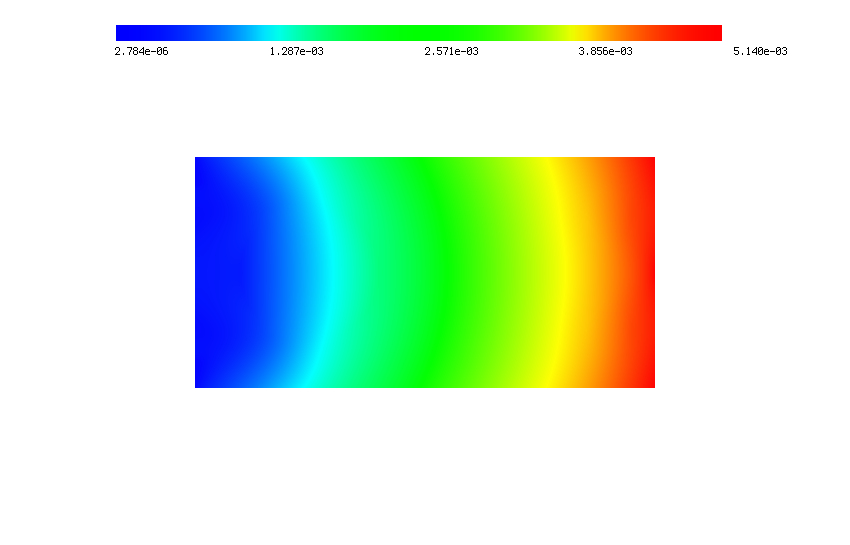} \\
        (a) & (b) \\
        \includegraphics[width=.5\textwidth]{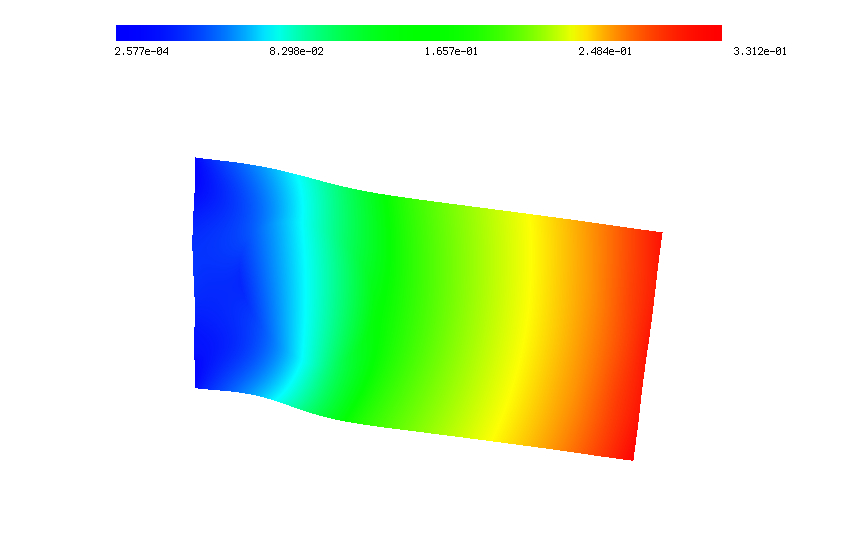} &
        \includegraphics[width=.5\textwidth]{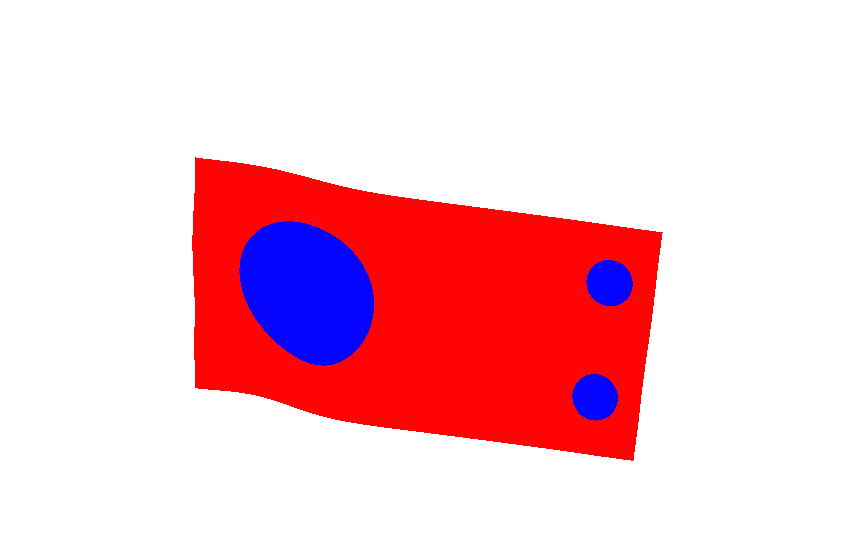} \\
        (c) & (d)
    \end{tabular}
    \caption{(a) Computational domain $\Dsf$ with subdomain $\Omega$ for Examples 4 and 5. (b) Solution to linear state equation of Example 4 (undeformed). (c) Solution to nonlinear state equation of Example 5 in deformed configuration. (d) Deformed domain of Example 5.}
    \label{fig_ela}
  \end{figure}

In this section, we consider vector-valued partial differential equations coming from elasticity as PDE constraints.
We consider the computational domain $\Dsf = (0,2)\times (0,1)$ which is clamped at the top left and bottom left, $\Gamma_D = \{(x_1, x_2): x_1 = 0 \mbox{ and } x_2 \in [0, 0.12] \cup [0.88,1]\}$ and is subject to a downward directed force at $\Gamma_N = {1}\times (0.45, 0.55)$. On the rest of the boundary, $\Gamma :=  \partial \Dsf \setminus (\Gamma_D \cup \Gamma_N)$, homogeneous Neumann boundary conditions are set. The domain consists of a strong material with Young's modulus $E = 1000$, and a very weak material which mimicks void inside $\Omega = B((0.5, 0.5)^\top, 0.3) \cup B((1.8, 0.25)^\top, 0.1)\cup B((1.8, 0.75)^\top, 0.1)$ where $E =0.1$, see Figure \ref{fig_ela}(a).
The Poisson ratio is set to $\nu = \frac13$ in both subdomains.
We consider the problem to minimize the compliance of the structure which is subject to self-weight as well as an external force acting on the boundary $\Gamma_N$. The problem reads
 \begin{align} \label{eq_exEla_J}
    \underset{\Omega}{\mbox{min }} J(u, \Omega):= \frac12 \int_\Dsf S_\Omega(\nnabla u):\nnabla u \; \mbox dx &
\end{align}
such that 
\begin{align}
     u \in H^1_{\Gamma_D}(\Dsf)^d: \int_\Dsf S_\Omega(\nnabla u) : \nnabla \psi \; \mbox dx =& \int_\Dsf f_\Omega \cdot \psi \; \mbox dx + \int_{\Gamma_N} g_N \cdot \psi \; \mbox dS   \label{eq_exEla_pde}
  \end{align}
  for all $\psi\in H^1_{\Gamma_D}(\Dsf)^d$,  where $S_\Omega : \VR^{d \times d} \rightarrow \VR^{d \times d}$ represents a stress tensor and $f_\Omega$ is a piecewise constant vector function. Also this problem fits into the framework considered in Section \ref{sec_derivTD} by setting  
  \begin{align*}
    j^\Omega(x, u(x), \nnabla u(x)) =& \frac12 S_{\Omega}(\nnabla u(x)) : \nnabla u(x), \\
    A_1^\Omega(x, u(x), \nnabla u(x)) =& 0,\\
    A_2^\Omega(x, u(x), \nnabla u(x)) =& S_{\Omega}(\nnabla u(x)),\\
    F_1^\Omega(x) =& f_{\Omega}(x),\\
    F_2^\Omega(x) =& 0.
  \end{align*}
  
\subsubsection{Example 4: Linear elasticity in 2D}
We consider problem \eqref{eq_exEla_J}--\eqref{eq_exEla_pde} for the two-dimensional cantilever example introduced above under a linear stress-strain relation,
\begin{align}
    S_\Omega(\nnabla u) = 2 \mu_\Omega e(u) + \lambda_\Omega Tr(e(u)) I_2
\end{align}
with the piecewise defined Lam\'e coefficients 
\begin{align*}
    \mu_\Omega(x) = \chi_\Omega(x) \mu_1 + \chi_{\Dsf \setminus \Omega}(x) \mu_2, \qquad
    \lambda_\Omega(x) = \chi_\Omega(x) \lambda_1 + \chi_{\Dsf \setminus \Omega}(x) \lambda_2 ,
\end{align*}
corresponding to Young's modulus $E_2 = 1000$, $E_1 = 0.1$ and Poisson ratio $\nu_2 = \nu_1 = \frac{1}{3}$, i.e.,
\begin{align}\label{eq_defLame}
    \mu_i = E_i \frac{1}{2 (1+ \nu_i)} \qquad \lambda_i = E_i  \frac{\nu_i}{ (1+\nu_i)(1-2\nu_i)}, \quad \mbox{for }i = 1,2.
\end{align}
Here, $e(u) = \frac12 (\nnabla u + \nnabla u^\top)$ denotes the linearized strain tensor, $Tr(A)$ denotes the trace of a matrix $A$ and $I_2$ denotes the two-dimensional identity matrix.
Moreover, in this example we neglect self-weight, $f_\Omega = 0$, and consider the external load $g_N = (0,-1)^\top$ acting on $\Gamma_N$.

We computed the unperturbed state and adjoint on a mesh with about 150 000 vertices (resulting in about 300 000 degrees of freedom) which is highly refined around the point $z = (1.2, 0.5)^\top$ for which we compute and verify the topological derivative. We solve the corrector equation \eqref{eq_K} with four considered two-dimensional inclusion shapes $\omega_{2D}^{(i)}$, $i=1, \dots 4$ where we replace the unbounded domain $\VR^d$ by the large, but bounded domain $B_R = B(\mathbf 0, R)$ with $R = 1000$ using a mesh with about 370 000 vertices.
We again compute the quantities $\delta \mathcal J$ defined in \eqref{eq_deltaJ} for a range of values for the inclusion radius $\eps$ \eqref{eq_defEpsUnderline} with $\eps_0 = 0.005$ and $\delta = 1.5$. The behavior of $\delta \mathcal J$ for these five inclusion shapes can be seen in Figure \ref{fig_deltaJ_ex4} where we can again observe the behavior $\delta \mathcal J =  \mathcal O(\eps^3)$ for all five inclusion shapes.

\subsubsection{Example 5: Nonlinear elasticity in 2D}
Finally, we consider problem \eqref{eq_exEla_J}--\eqref{eq_exEla_pde} with the nonlinear St. Venant-Kirchhoff material
\begin{align}
    S_\Omega(\nnabla u) := (I_2 + \nnabla u) \left[ \lambda_\Omega \mbox{Tr}(\frac{1}{2}(C( \nnabla u)-I_2))I_2 + \mu_\Omega(C(\nnabla u)-I_2) \right],
\end{align}
with $C(\nnabla u ):= (I_2+\nnabla u)^\top (I_2+\nnabla u )= I_2+ \nnabla u  + \nnabla u^\top + \nnabla u^\top \nnabla u$, and $\lambda_\Omega = \chi_\Omega \lambda_1 + \chi_{\Dsf \setminus \Omega} \lambda_2$, $\mu_\Omega= \chi_\Omega \mu_1 + \chi_{\Dsf \setminus \Omega} \mu_2$ with $E_2 = 1000$, $E_1 = 0.1$ and $\nu_2 = \nu_1 = 0.3$ and $\lambda_1, \lambda_2, \mu_1, \mu_2$ defined as in \eqref{eq_defLame}. Here we use the larger external force $g_N = (0, -20)^\top $ acting on $\Gamma_N$ and also mimic self-weight with $f_2 = (0,-5)^\top$, $f_1 = 0$.
  
For solving the state, adjoint and corrector equations, we use the same meshes as in Example 4. For solving the nonlinear state equation \eqref{eq_abstr_pde} we use a load stepping scheme with $N_{ls} = 20$ load steps where we gradually increase the loads given by $f_2$ and $g_N$ until they reach their full values given above, i.e., we solve \eqref{eq_abstr_pde} with the data given above and
\begin{align*}
    f_2 =  \frac{k}{N_{ls}} \left( \begin{array}{c} 0 \\ -5 \end{array} \right) \qquad g_N = \frac{k}{N_{ls}} \left( \begin{array}{c} 0 \\ -20 \end{array} \right)
\end{align*}
for $k = 1, \dots N_{ls}$. For each load step, we use Newton's method to solve \eqref{eq_abstr_pde} for the given loads. The solution to the nonlinear state equation in deformed configuration can be seen in Figure \ref{fig_ela}(c) and the deformed domain in Figure \ref{fig_ela}(d). Moreover, we solve the corrector equation \eqref{eq_K} by means of a damped Newton method with a very conservative damping factor of $0.002$.

\begin{figure}
    \begin{center}\includegraphics[width=.5\textwidth]{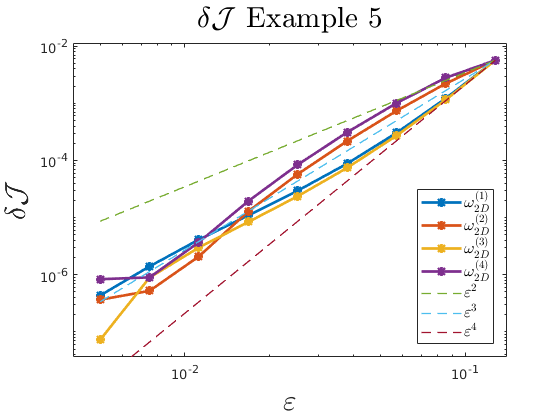}\end{center}
    \caption{Numerical verification of topological asymptotic expansion \eqref{eq_topAsympExp} for Example 5 for four different inclusion shapes.}
    \label{fig_deltaJ_ex5}
\end{figure}

We performed a Taylor test by computing the quantities $\delta \mathcal J$ for the four inclusion shapes $\omega_{2D}^{(i)}$, $i=1,2,3,4$, using the same parameters $\underline \eps$ as in Example 4. Figure \ref{fig_deltaJ_ex5} shows that $\delta \mathcal J$ decays at least as fast as $\eps^3$ which confirms the topologcial asymptotic expansion \eqref{eq_topAsympExp}.
  
\section{On the efficient evaluation of topological derivatives in the full domain} \label{sec_effEval}
In Section \ref{sec_numExp}, we applied the systematic procedure for the computation of topological derivatives introduced in Section \ref{sec_derivTD} to different model problems involving linear and nonlinear PDE constraints in two and three space dimensions. Moreover we verified the obtained values by Taylor tests. However, looking at the formulas of Section \ref{sec_summary}, it can be seen that the evaluation of the topological derivative at a spatial point $z$ requires the solution of the corrector equation \eqref{eq_K}. On the other hand, in order to employ an optimization algorithm, one usually is interested in the topological derivative in the full design domain. Of course, solving \eqref{eq_K} for every spatial point in the design domain (or every element of a mesh) is prohibitively expensive. 

In the following, we show how, for linear or semilinear PDE constraints, the corresponding solution $K$ can be obtained as a linear combination of some functions which can be precomputed. Note that, since we assume the PDE constraint to be linear or semilinear, both $A_2^{\inn}$ and $A_2^{out}$ are linear with respect to the third argument $\nnabla u$ and constant with respect to the second argument $u$, i.e., there exist $a_2^{\inn}(z), a_2^{out}(z)$ such that
\begin{align*}
    A_2^{\inn}(z, u_0(z), \nnabla u_0(z)) = a_2^{\inn}(z) \nnabla u_0(z) \quad \mbox{and} \quad A_2^{out}(z, u_0(z), \nnabla u_0(z)) = a_2^{out}(z) \nnabla u_0(z).
\end{align*}
We introduce the notation
\begin{align*}
    a_2^\omega(z, \nnabla u_0(z)) := \chi_\omega(x) a_2^{\inn}(z) \nnabla u_0(z) + \chi_{\VR^d \setminus \omega}(x) a_2^{out}(z) \nnabla u_0(z).
\end{align*}
Now let $\hat K$ be the solution to
\begin{align} \label{eq_defKHat}
    \int_{\VR^{d}} &   a_2^\omega(z,\nnabla \hat K)  : \nnabla \psi \;dx
    = \int_{\omega}  (F_2^{\inn}-F_2^{out})(z) : \nnabla \psi  \;dx 
\end{align}
for all $\psi$ and, for $i \in \{1,\dots, m\}$, $j \in \{1, \dots, d\}$,  let $\tilde K_{e_{ij}}$ be the solution to
\begin{align} \label{eq_defKTilde}
\begin{aligned}
    \int_{\VR^{d}} &   a_2^\omega(z, \nnabla \tilde K_{e_{ij}}) : \nnabla \psi \;dx
                   &=  -  \int_{\omega} ( a_2^{\inn} - a_2^{out})(z) e_{ij}  : \nnabla \psi \;dx
\end{aligned}
\end{align}
for all $\psi$. Here $e_{ij} \in \VR^{m \times d}$ denotes the unit basis element with the value $1$ at the position $(i,j)$ and the value $0$ else.

These considerations immediately yield the following lemma.

\begin{lemma} \label{lem_Kprecomp}
    Assume that the PDE constraint \eqref{eq_abstr_pde} is linear or semilinear and let $\nnabla u_0(z) \in \VR^{m \times d}$ be given. Then the solution $K$ to \eqref{eq_K} is given as the linear combination
    \begin{align}
        K =  \hat K + \sum_{i=1}^m \sum_{j=1}^d (a_2^{\inn}-a_2^{out})(z) \nnabla u_0(z)[i,j] \, \tilde K_{e_{ij}} 
    \end{align}
    with $\hat K$ and $\tilde K_{e_{ij}}$ given in \eqref{eq_defKHat} and \eqref{eq_defKTilde}, respectively.
\end{lemma}

Using Lemma \ref{lem_Kprecomp} it is possible to compute the solution $K$ of \eqref{eq_K} for any value of $\nnabla u_0(z)$ by linear combination once the $md+1$ functions $\hat K$, $\tilde K_{e_{ij}}$ have been computed. By integration over $\omega$ or $\VR^d$ (which is numerically approximated by $B(\mathbf 0, R)$ with $R=1000$), the topological derivative can then be approximately evaluated at every point without having to solve an additional boundary value problem.

\begin{figure}
    \begin{tabular}{cc}
        \includegraphics[width=.5\textwidth]{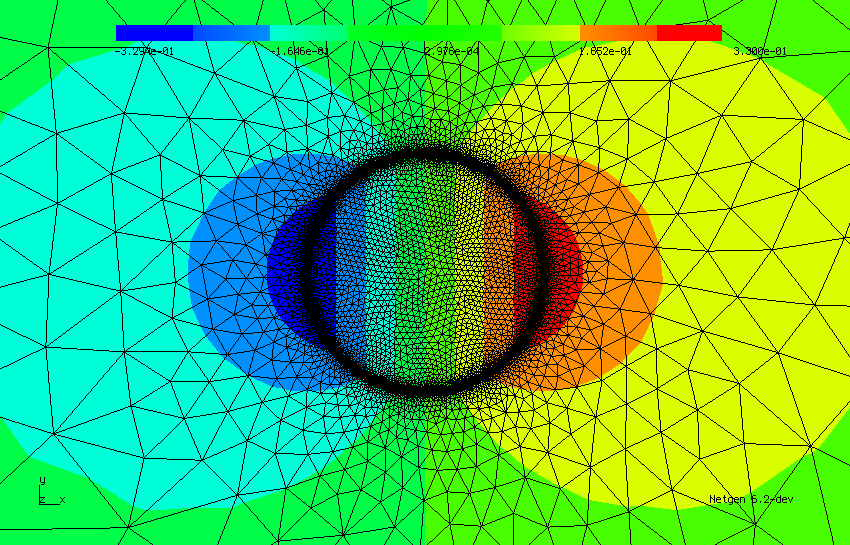}&
        \includegraphics[width=.5\textwidth]{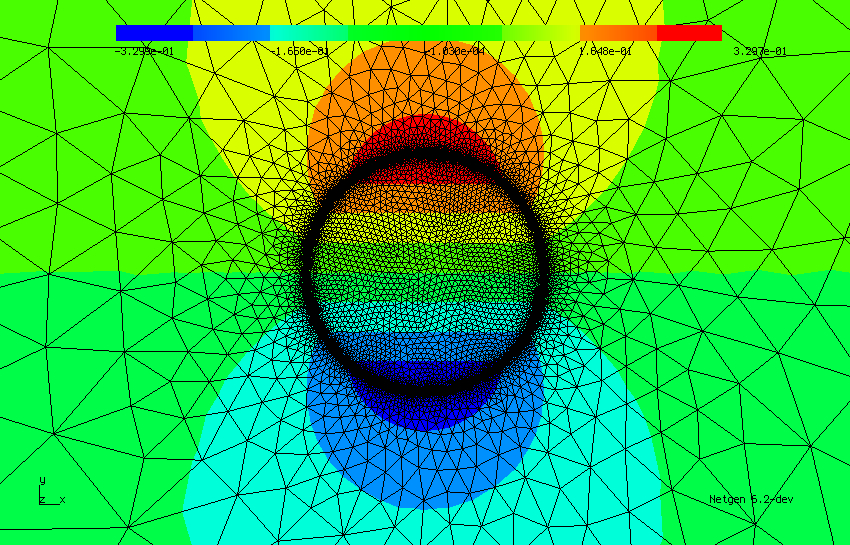}
    \end{tabular}
    \caption{Solutions to \eqref{eq_defKTilde} for $e_{11} = (1,0)^\top$ (left) and $e_{12} = (0,1)^\top$ (right).}
    \label{fig_Ktilde}
\end{figure}

In the case of Example 1 (Section \ref{sec_2dlin}), it follows from $\mathbf M_1 = \mathbf M_2 = (0,0)^\top$ that $\hat K = 0$. Since $m=1$, $d=2$, we need to precompute the solutions to \eqref{eq_defKTilde} for the two unit vectors in $\VR^2$ and can obtain the solution $K$ to \eqref{eq_K} by Lemma \ref{lem_Kprecomp}. The solutions $\tilde K_{e_{11}}$ and $\tilde K_{e_{12}}$ to \eqref{eq_defKTilde} with $\VR^d$ replaced by $B(\mathbf 0, 1000)$ are depicted in Figure \ref{fig_Ktilde}.

\begin{figure}
    \begin{tabular}{ccc}
        \includegraphics[width=.33\textwidth, trim = 70 0 60 0, clip]{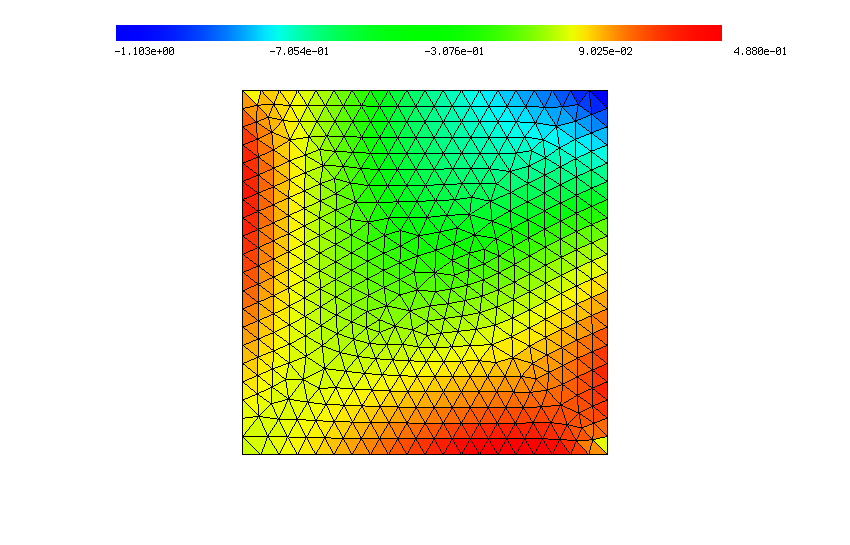} & 
        \includegraphics[width=.33\textwidth, trim = 70 0 60 0, clip]{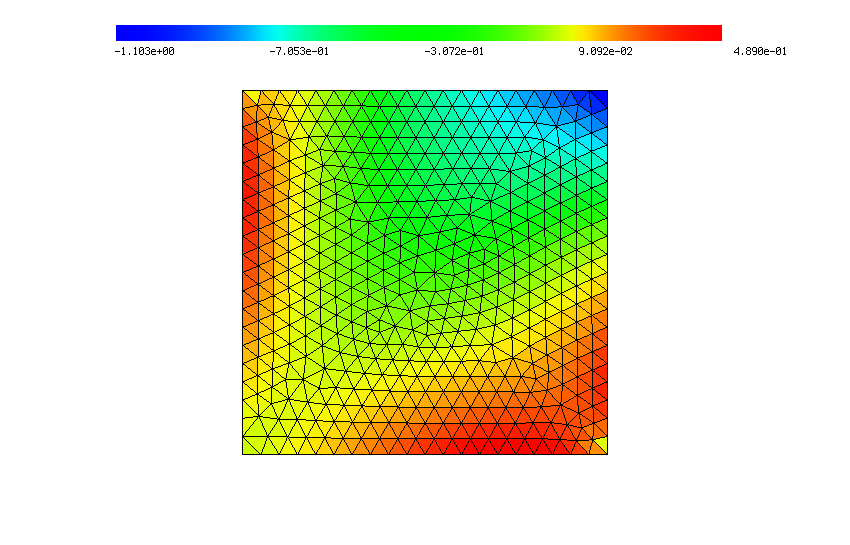} & 
        \includegraphics[width=.33\textwidth, trim = 70 0 60 0, clip]{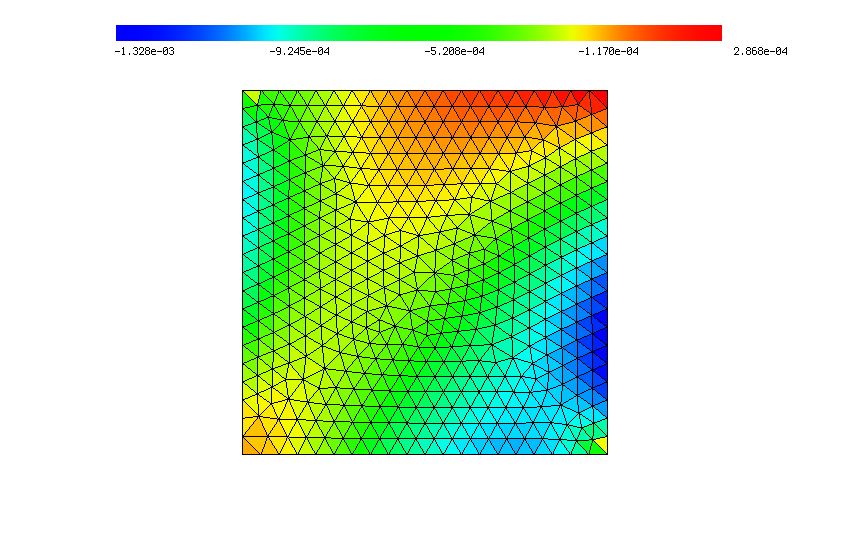} 
    \end{tabular}
    \caption{Left: Numerically computed topological derivative for Example 1. Center: Analytical formula \eqref{eq_TD_Ex1_closedForm}. Right: Difference between numerical and analytical formulas. Observe the different orders of magnitude.}
    \label{fig_ex1_numVSana}
\end{figure}

Next, we chose a (comparably coarse) mesh with only background material, i.e. $\Omega = \emptyset$, and numerically computed the topological derivative for the centroid of every triangular element by evaluating the formulas of Section \ref{sec_summary}. Also here, integrals over $\VR^d$ were approximated by integrals over the large, but bounded ball $B(\mathbf 0, 1000)$. A comparison of the obtained results with the analytical formula for the topological derivative given in \eqref{eq_TD_Ex1_closedForm} showed good accordance, see Figure \ref{fig_ex1_numVSana}.

\begin{remark}
The procedure described in this section is still computationally expensive since it involves the numerical computation of several integrals over $\omega$ and $B(\mathbf 0, 1000)$ for each topological derivative evaluation.

We remark that in many cases it is possible to exploit the (affine) linearity of the topological derivative $d \mathcal J$ \eqref{eq_summary_dJ} with respect to the adjoint state as well as some rotational symmetry property of $K$ with respect to $\nnabla u_0(z)$ to directly precompute the topological derivative $d \mathcal J$ for some basis elements and to obtain the full topological derivative as a linear combination of these values. We remark that such a procedure -- if applicable -- can be computationally cheaper and allow for the use of the obtained formulas in iterative topology optimization algorithms, see e.g. \cite{a_AMGA_2019a, GanglSturm2021Hcurl} for more details.
\end{remark}

\section*{Conclusion and Outlook}
We have illustrated a systematic way of computing topological derivatives for a large class of PDE-constrained topology optimization problems. Using automated differentiation in \texttt{NGSolve}, it was possible to compute all potential terms of the topological derivative for an abstract problem class. We showed the effectivity of the method in several linear and nonlinear model problems and verified the topological derivative formulas by computing the values of the topological asymptotic expansion that should be satisfied by the topological derivative.

The presented work can be extended in several directions. On the one hand, we only considered problems posed in an $H^1$ setting. We remark that the same procedure is possible for any differential operator $\nnabla$ that scales as $(\nnabla \varphi)\circ T_\eps = \frac{1}{\eps} \nnabla (\varphi \circ T_\eps)$. In particular, this also includes the $curl$ operator, thus also allowing to treat $H(curl)$ problems as they arise in electromagnetics \cite{GanglSturm2021Hcurl}. Moreover, for some linear problems, the computational effort for solving the corrector equation \eqref{eq_K} could be reduced by employing a boundary element method rather than a finite element method on an approximation of the unbounded domain. Finally, a further topic of future research is the numerical analysis of the topological derivative with respect to mesh sizes in the bounded domain $\Dsf$ and in the blown-up domain $B_R$ as well as with respect to the radius $R$ of the large domain.

\bibliographystyle{plain}
\bibliography{topoNLEla}

\end{document}